\documentclass[journal]{IEEEtran}
\usepackage{epstopdf}
\usepackage[cmex10]{amsmath}
\usepackage{amsfonts}
\usepackage{amssymb}
\usepackage{graphicx}
\usepackage{algorithm}
\usepackage{verbatim}
\usepackage{array}
\usepackage{multirow}
\usepackage{dcolumn}
\usepackage{color}
\usepackage[noadjust]{cite}
\usepackage{url}
\usepackage{balance}
\usepackage{booktabs}
\DeclareGraphicsExtensions{.eps}
\usepackage{setspace}
\usepackage{algpseudocode}
\usepackage[export]{adjustbox}
\usepackage{adjustbox}
\usepackage{soul}

\def\u{\mathbf{u}}

\def\U{\mathcal{U}}

\algnewcommand{\IIf}[1]{\State\algorithmicif\ #1\ \algorithmicthen}
\algnewcommand{\EndIIf}{\unskip\ \algorithmicend\ \algorithmicif}

\begin{document}





\title{Two-Stage Robust Unit Commitment for Co-Optimized Electricity Markets: An Adaptive Data-Driven Approach for Scenario-Based Uncertainty Sets}




\author{Alexandre~Velloso,~\IEEEmembership{Student Member,~IEEE,}
		Alexandre~Street,~\IEEEmembership{Senior Member,~IEEE,}
		David~Pozo,~\IEEEmembership{Member,~IEEE,}	
	    Jos{\'e}~M.~Arroyo,~\IEEEmembership{Senior Member,~IEEE,}	
		Noemi~G.~Cobos,~\IEEEmembership{Student Member,~IEEE}	
		

\thanks{The work of Alexandre Velloso was supported by Finep/PIPG -- Financiadora de Estudos e Projetos / Programa de Incentivo \`{a} P\'{o}s-Gradua\c{c}\~{a}o. Alexandre Street would like to acknowledge the financial support from CNPq -- Conselho Nacional de Desenvolvimento Cient\'{i}fico e Tecnol\'{o}gico -- and FAPERJ -- Funda\c{c}\~{a}o de Amparo \`{a} Pesquisa do Estado do Rio de Janeiro. Jos\'{e} M. Arroyo would like to acknowledge the financial support from the Ministry of Science, Innovation and Universities of Spain, under Project RTI2018-098703-B-I00 (MCIU/AEI/FEDER, UE). The work of D. Pozo was supported by Skoltech NGP Program (Skoltech-MIT joint project)}
}

\maketitle


\begin{abstract}

	Two-stage robust unit commitment (RUC) models have been widely used for day-ahead energy and reserve scheduling under high renewable integration. The current state of the art relies on budget-constrained polyhedral uncertainty sets to control the conservativeness of the solutions. The associated lack of interpretability and parameter specification procedures, as well as the high computational burden exhibited by available exact solution techniques call for new approaches. In this work, we use an alternative scenario-based framework whereby uncertain renewable generation is characterized by a polyhedral uncertainty set relying on the direct specification of its vertexes. Moreover, we present a simple, yet efficient, adaptive data-driven procedure to dynamically update the uncertainty set vertexes with observed daily renewable-output profiles. Within this setting, the proposed data-driven RUC ensures protection against the convex hull of realistic scenarios empirically capturing the complex and time-varying intra-day spatial and temporal interdependences among renewable units. The resulting counterpart features advantageous properties from a computational perspective and can be effectively solved by the column-and-constraint generation algorithm until $\epsilon$-global optimality. Out-of-sample experiments reveal that the proposed approach is capable of attaining efficient solutions in terms of cost and robustness while keeping the model tractable and scalable.

\end{abstract}



\section*{Nomenclature}

This section lists the notation. Bold symbols are reserved to matrices (uppercase) and vectors (lowercase). Additional symbols with superscript ``$(k)$" denote new variables corresponding to the $k$-th scenario selected by the solution method. 

\subsection{Sets}

\begin{description}
	[\setlabelwidth{$\mbox{\it{ELLLLL}}$}\usemathlabelsep]
	\vspace{0.02cm}
	\item[$\mathcal{F}$] Feasible set for the decision variables associated with thermal generators.
	
	\vspace{0.02cm}
	
	\item[$\mathcal{H}$] Set of time periods.  
	
	\vspace{0.02cm}

	\item[$\mathcal{R}$] Set of renewable energy generators.
	\vspace{0.02cm}
	\item[$\mathcal{S}$] Set of renewable energy generation scenarios.
	\vspace{0.02cm}
	\item [$\mathcal{S}_j$] Subset of $\mathcal{S}$ including the scenarios selected until iteration $j$.
	
	\vspace{0.02cm}
	\item[$\U$] Uncertainty set defined as the convex hull of $\mathcal{S}$.
	
	\vspace{0.02cm}
	
	\item[$\mathcal{X}$] Set of first-stage generation-related variables.
	
\end{description}

\subsection{Constants}

\begin{description}
	[\setlabelwidth{$\mbox{\it{ELLLLL}}$}\usemathlabelsep]
	
	\vspace{0.02cm}
	\item[$\Gamma,\,\Lambda$] Spatial and temporal budgets used in the budget-constrained uncertainty set formulation.
	\vspace{0.02cm}
	\item[$\Delta_{ih}^{+}, \Delta_{ih}^{-}$] Upper and lower deviation bounds for renewable generator $i$ and period $h$ used in the budget-constrained uncertainty set formulation.%
	\vspace{0.02cm}

	\item[$\epsilon$] Feasibility tolerance.
	\vspace{0.02cm}
	\item[$\mathbf{A}$] Line-bus incidence matrix.
	\vspace{0.02cm}
	\item[$\mathbf{B}$] Thermal generator-bus incidence matrix.
	
    \vspace{0.02cm}
	\item[$\mathbf{c}^{g}$] Vector of fuel costs of thermal generators.
	\vspace{0.02cm}
	\item[$\mathbf{c}^{dn},\mathbf{c}^{up}$] Vectors of cost rates  for down- and up-spinning reserves.
    \vspace{0.02cm}
	\item[$\mathbf{{d}}_{h}$] Vector of nodal consumptions in period $h$.
	\vspace{0.02cm}
	\item[$\mathbf{e}$] Vector of ones with appropriate dimension.
	\vspace{0.02cm}
	\item[$\mathbf{\overline{f}}$] Vector of line capacities.
	\vspace{0.02cm}
	\item[$\mathbf{\overline{G}},\mathbf{\underline{G}}$] Diagonal matrices of maximum and minimum generation limits of thermal units.
	\vspace{0.02cm}
	\item[$\overline{G}_i$] Capacity of thermal unit $i$.
	\vspace{0.02cm}
	\item[$K$] Number of days for observed renewable energy generation data.
	\vspace{0.02cm}
	\item[$\mathbf{P}$] Renewable generator-bus incidence matrix.
	\vspace{0.02cm}
	
	\item[$RD_i,RU_i$] Ramp-down and ramp-up limits of thermal unit $i$ within two consecutive periods.
	\vspace{0.02cm}
	\item[$\mathbf{R}^{dn},\mathbf{R}^{up}$] Diagonal matrices of downward and upward limits for corrective actions of thermal units within each period.
	\vspace{0.02cm}
	\item[$\mathbf{S}$] Angle-to-flow matrix.
	\vspace{0.02cm}
	\item[$SD_i,SU_i$] Shut-down and start-up ramp rates for thermal unit $i$.
	\vspace{0.02cm}
	
	\item[$\mathbf{u}_{hk}$] Vector of renewable energy generation in period $h$ for scenario $k$. 	
	
	\vspace{0.02cm}
	\item[$\mathbf{U}_k$] Scenario $k$ of $\mathcal{S}$.
	\vspace{0.02cm}
	\item[$\mathbf{\hat{u}}_{h}$] Vector of expected day-ahead renewable energy generation levels in period $h$.
	\vspace{0.02cm}
	\item[$\hat{u}_{ih}$] Expected day-ahead power output of renewable generator $i$ in period $h$.
	\vspace{0.02cm}

\end{description}

\subsection{First-Level Decision Vectors}
\begin{description}[\setlabelwidth{$\mbox{\it{ELLLLL}}$}\usemathlabelsep]
	\vspace{0.02cm}
	\item[$\mathbf{\boldsymbol{\theta}}_h$] Phase angles in period $h$.
	\vspace{0.02cm}
	\item[$\mathbf{c}_{h}^{sd}, \mathbf{c}_{h}^{su}$] Shut-down and start-up costs in period $h$.
	\vspace{0.02cm}
	\item[$\mathbf{f}_{h}$] Line power flows in period $h$.
	\vspace{0.02cm}
	\item[$\mathbf{g}_{h}$] Nominal generation levels in period $h$.
	\vspace{0.02cm}
	\item[$\mathbf{r}^{dn}_{h},\mathbf{r}^{up}_{h}$] Down- and up-spinning reserves in period $h$.
	\vspace{0.02cm}
	\item[$\mathbf{v}_{h}$] Binary on/off statuses in period $h$.
	

\end{description}

\subsection{Second- and Third-Level Decision Variables}
\begin{description}[\setlabelwidth{$\mbox{\it{ELLLLL}}$}\usemathlabelsep]
	\vspace{0.02cm}
	\item[$\alpha_k$] Second-level decision variable representing the convex combination weight of scenario $\mathbf{U}_k$ used in the scenario-based uncertainty set $\mathcal{U}$.
	\vspace{0.02cm}
	\item[$\mathbf{\boldsymbol{\theta}}_h^{wc}$] Third-level  decision vector representing the worst-case nodal phase angles in period $h$.
	\vspace{0.02cm}
	\item[$\mathbf{f}_{h}^{wc}$] Third-level decision vector representing the worst-case line power flows in period $h$.
	\vspace{0.02cm}
	\item[$\mathbf{g}_{h}^{wc}$] Third-level decision vector representing the worst-case generation redispatch in period $h$.
	\vspace{0.02cm}
	\item[$\mathbf{s}^{+}_h,\mathbf{s}^{-}_h$] Third-level decision vectors representing  the renewable energy spillage and load shedding in period $h$.
	\vspace{0.02cm}
	\item[$\mathbf{U}$] Second-level decision matrix representing the generation levels for all renewable units across the day-ahead scheduling horizon.
	\vspace{0.02cm}
	\item[$\mathbf{u}_h$] Second-level decision vector representing the $h$-th column of $\mathbf{U}$.
	\vspace{0.02cm}
	\item[$u_{ih}$] Second-level decision variable representing the $i$-th element of $\mathbf{u}_h$.
	\vspace{0.02cm}
	\item[$z_{ih}^{+}, z_{ih}^{-}$] Second-level decision variables, used in the budget-constrained uncertainty set formulation, representing upward and downward deviations from the expected value for renewable generator $i$ in period $h$.

\end{description}

\subsection{Dual Variables}

\begin{description}[\setlabelwidth{$\mbox{\it{ELLLLL}}$}\usemathlabelsep]
	\vspace{0.02cm}
	\item[$\boldsymbol{\beta}_{h}$] Vector of dual variables associated with the power balance equations in period $h$.
	\vspace{0.02cm}
	\item[$\boldsymbol{\gamma}_h,\boldsymbol{\tau}_h$] Vectors of dual variables associated with generation redispatch limits in period $h$.
	\vspace{0.02cm}
	\item[$\boldsymbol{\zeta}_h, \boldsymbol{\kappa}_h$] Vectors of dual variables associated with ramping-up and ramping-down constraints in period $h$.
	\vspace{0.02cm}
	\item[$\boldsymbol{\sigma}_h,\boldsymbol{\pi}_h$] Vectors of dual variables associated with the power flow capacity constraints in period $h$.
	
	
	\vspace{0.02cm}
	\item[$\boldsymbol{\varsigma}_h, \boldsymbol{\xi}_h$] Vectors of dual variables associated with load shedding and renewable energy spillage limits in period $h$.
	\vspace{0.02cm}
	\item[$\boldsymbol{\omega}_h$] Vector of dual variables associated with Kirchhoff's second law constraints in period $h$.

\end{description}

\subsection{Functions}

\begin{description}
	[\setlabelwidth{$\mbox{\it{ELLLLL}}$}\usemathlabelsep]
	\vspace{0.02cm}
	\item[$\Phi(\cdot)$] Worst-case system power imbalance.
	\vspace{0.02cm}
	\item[$\mathbf{a(\cdot)} , \mathbf{b(\cdot)} $] Vector functions used to enforce ramping limits.
	
	
	
\end{description}

\vspace{0.02cm}

\section{Introduction}\label{Introduction}


\PARstart{N}{on-dispatchable} renewable energy generation (REG) has undergone a sharp increase in the last decades and is already one of the major components in some electricity markets. 
High integration of these intermittent and variable energy sources brings additional challenges to short-term power system operation that are well known and have been widely discussed \cite{zheng2015stochastic,Jiang2012,bertsimas2013adaptive,yongpeiEJOR2014,lorca2014adaptive,moreira2015energy,BoZeng2015,dai2016multi,wang2017robust,shahidehpour2017security,Noemi2017,cobos2018network}. Briefly, REG variability, especially from wind power units, is driven by complex time-varying spatial and temporal dynamics \cite{giebel2011state}. In order to benefit from REG resources, a mix of conventional generation and expensive operational actions are both needed to constantly deploy (up and down) reserves in a fast and reliable way. For this reason, the uncertainty inherent to REG should be precisely accounted for in the scheduling and dispatch models used to determine appropriate levels of energy and reserves, such as those adopted in currently implemented co-optimized electricity markets \cite{moreira2015energy,arroyo2005energy,wang2009contingency,chen2014incorporating,Noemi2017,cobos2018network,FERC, GREECE}.

Due to the appealing tradeoff between tractability and accuracy, two-stage adaptive robust optimization has been used to deal with uncertainty in day-ahead generation scheduling  \cite{Jiang2012,bertsimas2013adaptive,yongpeiEJOR2014,moreira2015energy,BoZeng2015,dai2016multi,shahidehpour2017security,wang2017robust,Noemi2017,lorca2014adaptive,cobos2018network}. The interested reader is referred to \cite{ShabbirBook} for a detailed literature review. In such robust unit commitment (RUC) models, a trilevel optimization problem is built to characterize the min(decision)-max(uncertainty)-min(decision) structure. Within such a scheme, the first-level problem determines, before the observation of the uncertain parameters, the day-ahead commitment for each generator. In the second-level problem, the worst-case scenario of uncertainties is selected within a given polyhedral uncertainty set as a function of the first-level decisions. Finally, in the third level, the best operational reaction (redispatch) is obtained for the second-level scenario within the first-level scheduled resources.

A key aspect in RUC is the way that uncertainties are represented. The success of an RUC model mainly depends on the selection of an uncertainty set that is capable of capturing the main patterns present in the uncertain parameters while keeping model tractability. Valuable and thorough discussions on the subject can be found in \cite{ShabbirBook} and \cite{guan2014uncertainty}.

Most previously reported two-stage RUC models  \cite{Jiang2012,bertsimas2013adaptive,yongpeiEJOR2014,moreira2015energy,BoZeng2015, dai2016multi,shahidehpour2017security,wang2017robust,Noemi2017,lorca2014adaptive, cobos2018network} rely on the budget-constrained polyhedral uncertainty set presented in \cite{bertsimas2004price}. The specification of a budget-constrained uncertainty set is made through linear inequality constraints defining the boundaries of a polyhedron. Such boundaries are set up by componentwise box-like limits and linear (budget) constraints limiting the number of components deviating from their nominal scenario. In this framework, each vertex of the polyhedron representing the uncertainty set is indirectly determined by the intersection of constraints, which may hinder the physical interpretation of the scenarios.

Variants of conventional budget-constrained uncertainty sets have been proposed in \cite{Bandi2012} and applied in {\cite{lorca2014adaptive,moreira2015energy,BoZeng2015,dai2016multi,wang2017robust,shahidehpour2017security}} to better model the variability of renewable power generation.  In {\cite{lorca2014adaptive}}, an interesting approach based on linear models was proposed to improve the dynamics of REG scenarios. The use of parametric linear models was also described in {\cite{moreira2015energy}}. For those cases, the polyhedral uncertainty sets were defined based on affine constraints representing linear models for renewable injections. In addition, budget constraints were applied to the residuals of the models to control the conservativeness of solutions. In \cite{BoZeng2015}, An and Zeng introduced variants of robust models based on linear expected-value operators (averages) applied to multiple worst-case operational costs. In such a work, each worst case resulted from a different budget-constrained uncertainty set centered on an exogenously generated scenario. Multiple spatial and temporal budget constraints were presented in {\cite{dai2016multi}} to increase the modeling capability. In {\cite{wang2017robust}}, a flexible uncertainty set was characterized in terms of a user-defined parameter to capture the risk of misestimating the box-like limits for REG levels. In {\cite{shahidehpour2017security}}, a flexible uncertainty set was adjusted over time to provide a tradeoff between economics and robustness of the generation schedule. Notwithstanding, the approaches presented in \cite{lorca2014adaptive,moreira2015energy,BoZeng2015,dai2016multi,wang2017robust,shahidehpour2017security} rely on budget-constrained uncertainty sets. Therefore, the modeling choices for describing REG variability, which features complex, nonlinear, and time-varying dynamics \cite{giebel2011state}, are restricted to linear models due to tractability issues.

From a computational perspective, another drawback associated with the use of budget-constrained uncertainty sets is the combinatorial growth of the number of vertexes of the polyhedron with respect to the number of uncertain parameters. The resulting RUC models based on budget-constrained uncertainty sets are challenging instances of trilevel programming that are generally addressed through decomposition-based methods. The state-of-the-art techniques such as the column-and-constraint generation algorithm (CCGA) {\cite{BoZeng2011}} and Benders decomposition {\cite{bertsimas2013adaptive}} involve the iterative solution of a master problem and a subproblem, also known as oracle subproblem. In the related literature, the oracle subproblem represents the worst-case-scenario search procedure corresponding to the two lowermost optimization levels of the trilevel counterpart. Hence, the oracle is an instance of bilevel programming, which, in the presence of budget-constrained uncertainty sets, is NP-hard \cite{bertsimas2013adaptive}. Solution techniques available to tackle the oracle subproblem can be categorized in two groups. On the one hand, heuristic yet efficient methods, such as variants of the outer-approximation algorithm \cite{bazaraa2013nonlinear}, were applied in \cite{Jiang2012,bertsimas2013adaptive,yongpeiEJOR2014,lorca2014adaptive} and \cite{dai2016multi}. On the other hand, in \cite{moreira2015energy,wang2017robust,shahidehpour2017security,BoZeng2015,Noemi2017,cobos2018network}, well-known linearization procedures, relying on the binary representation of the vertexes of the polyhedron characterizing the uncertainty set, were used to obtain exact yet computationally expensive single-level equivalents based on mixed-integer linear programming (MILP). Thus, existing solution methodologies for RUC models with budget-constrained uncertainty sets either efficiently provide a solution without being able to acknowledge optimality or rely on exact MILP-based NP-hard models that are challenging to solve in practice.

Motivated by the above issues of existing works \cite{Jiang2012,bertsimas2013adaptive,yongpeiEJOR2014,lorca2014adaptive,moreira2015energy,BoZeng2015,dai2016multi,wang2017robust,shahidehpour2017security,Noemi2017,cobos2018network} and the wide availability of REG data, the objective of this paper is to propose an alternative to the use of budget-constrained uncertainty sets in RUC. In this work, we consider the RUC for co-optimized electricity markets, i.e., the  centralized robust joint scheduling of energy and reserves targeting total cost minimization. Here, as suggested in \cite{Bertsimas2011}, the uncertainty characterization is directly connected to data. To that end, we propose modeling the REG uncertainty in day-ahead RUC by an alternative scenario-based polyhedral uncertainty set that is built through a novel data-driven approach. Based on the general scenario-based uncertainty set description provided in \cite{Bertsimas2009}, we define a new polyhedral uncertainty set as the convex hull of a set of exogenously generated multivariate points, or scenarios, capturing relevant information regarding the uncertainty process over a given time window. Thus, differently from \cite{BoZeng2015}, each vertex of the polyhedron representing the uncertainty set is defined as one of these exogenous scenarios. In the proposed data-driven framework, scenarios represent observed daily renewable-generation profiles, i.e., matrices whose dimension is given by the number of renewable units and the number of time periods of the scheduling horizon, typically 24 hours. Hereinafter, the proposed data-driven scenario-based uncertainty set is referred to as DDUS.

Two recent examples of successful application of the scenario-based uncertainty sets first proposed in \cite{Bertsimas2009} can be found in \cite{Fernandes2016961} and \cite{Bertsimas2017}. Within a finance context, the vertexes of the uncertainty set were generated directly from most recent observed data in \cite{Fernandes2016961}. Using a general mathematical setting, in \cite{Bertsimas2017}, sampled points were endogenously selected to belong to the uncertainty set through embedded statistical hypothesis tests. In the context of RUC, however, this alternative framework has not been explored yet despite its relevant benefits. 

From a modeling perspective, the use of the proposed uncertainty characterization for RUC is advantageous in several aspects as compared with previous models relying on budget-constrained uncertainty sets  \cite{Jiang2012,bertsimas2013adaptive,yongpeiEJOR2014,lorca2014adaptive,moreira2015energy,BoZeng2015,dai2016multi,wang2017robust,shahidehpour2017security,Noemi2017,cobos2018network}. First, the true underlying uncertainty process drives the construction of the polyhedron representing the uncertainty set, which is made up of vertexes with high physical interpretability. As a consequence, the resulting RUC features relevant information about the complex and time-varying temporal and spatial dependences found in REG within the scheduling horizon. Moreover, the novel data-driven procedure devised to build polyhedral uncertainty sets through their vertexes is an entirely exogenous adaptive and nonparametric process. Therefore, the proposed approach paves the way for the use of a wide range of existing scenario-generation methods as alternatives to the proposed data-driven procedure. For instance, any nonlinear model or data-processing scheme useful for defining or preprocessing the scenarios, such as clustering, data categorization, or filtering processes based on weather-related and real-time dispatch information, can be used to generate vertexes for the uncertainty set.

The incorporation of the proposed DDUS in RUC is also beneficial from a methodological perspective. Similar to existing models \cite{Jiang2012,bertsimas2013adaptive,yongpeiEJOR2014,lorca2014adaptive,moreira2015energy,BoZeng2015,dai2016multi,wang2017robust,shahidehpour2017security,Noemi2017,cobos2018network}, the resulting data-driven formulation, denoted by DDRUC, is suitable for the state-of-the-art CCGA. In addition, as a salient feature, the proposed scenario-based robust framework is characterized by a relevant property: one of the multivariate points within the given time window of observed data is the worst-case vertex provided by the optimal solution of the oracle subproblem. This property results in an oracle that is solvable in polynomial time \cite{Bertsimas2017}, unlike the oracle subproblems described in \cite{Jiang2012,bertsimas2013adaptive,yongpeiEJOR2014,lorca2014adaptive,moreira2015energy,BoZeng2015,dai2016multi,wang2017robust,shahidehpour2017security,Noemi2017,cobos2018network}. It is also worth mentioning that, according to our empirical findings, which are consistent with those reported in \cite{Fernandes2016961}, the use of DDUS typically requires a narrow time window to attain high-quality solutions. This aspect is particularly relevant for the practical adoption by system operators to schedule generation in electricity markets.

The contributions of this paper are threefold:
\begin{enumerate}
	\item To raise awareness of the modeling capability and computational advantages of the scenario-based polyhedral uncertainty sets proposed in \cite{Bertsimas2009} to address REG uncertainty in RUC problems. Within this general framework, we propose defining scenarios as matrices representing the hourly generation profiles of all renewable units within a day. Hence, for the first time in the RUC literature, REG variability is described by a polyhedral uncertainty set relying on the convex hull of a polynomial set of multivariate points.
	
	\item To propose a nonparametric data-driven procedure that defines the vertexes of the resulting polyhedral uncertainty set directly from observed data, thereby embedding the true complex and time-varying interdependences among renewable units.
	
	\item To present a novel data-driven two-stage robust unit commitment model that is scalable, easy to specify, and suitable for the exact CCGA due to the resulting computationally inexpensive and polynomial-time-solvable oracle subproblem.
\end{enumerate}

The rest of this paper is organized as follows. In Section \ref{sec.UncertaintySetFormulation}, the alternative uncertainty set is described. The proposed RUC model is formulated in Section \ref{sec.Model}. The solution methodology is presented in Section \ref{sec.Meth}. An evaluation procedure for assessing the devised model and numerical experience are reported in Section \ref{sec.Case}. Finally, this paper is concluded in Section \ref{sec.Conclusion}.

\vspace{-0.20cm}

\section{Uncertainty Set Characterization}
\label{sec.UncertaintySetFormulation}

\vspace{-0.05cm}

In generation scheduling under high penetration of renewable-based generation, uncertain data comprise the day-ahead power output for each renewable generator $i \in \mathcal{R}$ and each time period $h \in \mathcal{H}$. Thus, we represent uncertainty by a renewable generator-by-time matrix, $\mathbf{U} \in \mathbb{R}^{ \;\mathcal{|R|} \times \mathcal{|H|}}$, where $\mathbf{U} = [\mathbf{u}_1, ..., \mathbf{u}_h, ..., \mathbf{u}_{\mathcal{|H|}}]$ and $\mathbf{u}_h \in \mathbb{R}^{\mathcal{|R|}}$ is the uncertainty vector whose components $u_{ih}$ correspond to the available generation of each renewable generator $i$ in period $h$. Note that generation levels of different REG units may present spatial, temporal, and cross-lagged dependences.

In robust optimization, uncertainty is modeled through uncertainty sets. Before presenting the proposed DDUS, a general formulation for conventional budget-constrained uncertainty sets is provided based on the above matrix characterization.

\vspace{-0.3cm}
\subsection{Conventional Budget-Constrained Uncertainty Sets}
\label{sub.Conventional}
In robust generation scheduling, uncertainty sets typically rely on fluctuation intervals representing the support of the uncertainties \cite{Jiang2012,bertsimas2013adaptive,yongpeiEJOR2014,Noemi2017,cobos2018network}, and spatial and/or temporal budget constraints modeling the conservativeness or risk aversion of the decision maker \cite{Jiang2012,bertsimas2013adaptive,yongpeiEJOR2014,dai2016multi,Noemi2017,cobos2018network}. Thus, conventional budget-constrained uncertainty sets can be defined as the set of matrices $\mathbf{U}$ such that:
\begin{IEEEeqnarray}{ll}
	u_{ih} = \hat{u}_{ih} + \Delta_{ih}^{+}z_{ih}^{+} - \Delta_{ih}^{-}z_{ih}^{-}   \;\,\quad \forall i \in \mathcal{R}, \,\; \forall h \in \mathcal{H}\;\;\; \label{BUSlin1}\\
	0\leq  z_{ih}^{+}, z_{ih}^{-} \leq 1\,\qquad\qquad\qquad\quad \forall i \in \mathcal{R}, \,\; \forall h \in \mathcal{H}\;\;\; \label{BUSlin4}\\
	\sum_{i \in \mathcal{R}} z_{ih}^{+} + z_{ih}^{-} \leq \Gamma|\mathcal{R}| \qquad\qquad\quad \forall h \in \mathcal{H}  \label{BUSlin2}\\
	\sum_{h \in \mathcal{H}} z_{ih}^{+} + z_{ih}^{-} \leq \Lambda \;\;\; \qquad\qquad\quad\,\, \forall i \in \mathcal{R}. \label{BUSlin3}
\end{IEEEeqnarray}

As per \eqref{BUSlin1} and \eqref{BUSlin4}, decision variables $z_{ih}^{+}$ and $z_{ih}^{-}$ determine the available generation for renewable generator $i$ in period $h$, represented by $u_{ih}$. The lower and upper box-like limits, $\hat{u}_{ih} - \Delta_{ih}^{-}$ and $\hat{u}_{ih} + \Delta_{ih}^{+}$, are, respectively, the minimum and maximum possible values for $u_{ih}$. The first budget constraint \eqref{BUSlin2} limits, for each period, the number of renewable generating units that may deviate from their expected value by parameter $\Gamma$, which can vary from $0$ to $100\%$. As in \cite{bertsimas2004price}, this constraint means that up to $\lfloor\Gamma |\mathcal{R}|\rfloor$ generators\footnote{In cases where $\Gamma |\mathcal{R}|$ is not an integer, an additional renewable generator is allowed to deviate to make up for the residual fraction, $\Gamma|\mathcal{R}|-\lfloor\Gamma |\mathcal{R}|\rfloor$.} can simultaneously change their production from their expected values $\hat{u}_{ih}$ to their upper or lower bounds in the same period. The second budget constraint \eqref{BUSlin3} limits the number of periods in which each renewable generator $i\in\mathcal{R}$ can deviate from its expected production. For the sake of simplicity, parameter $\Lambda$ represents an integer number of periods ranging between $0$ and $|\mathcal{H}|$.
\vspace{-0.20cm}
\subsection{Scenario-Based Uncertainty Sets Driven by Data}
\label{DDUSsection}

In order to build a flexible, easy to specify, and practical uncertainty set that allows scalability for the available exact solution techniques for RUC, we propose a scenario-based polyhedral uncertainty set \cite{Bertsimas2009,Fernandes2016961,Bertsimas2017} to model REG uncertainty. To that end, the uncertainty set is defined as the convex hull of a set of multivariate points representing scenarios, i.e., REG profiles. Among the various exogenous approaches that could be used to generate the scenarios, we propose the use of a data-driven scheme. Historical daily profiles are therefore directly used as scenarios, thereby embedding relevant information about the true underlying uncertainty process in each vertex of the uncertainty set. Mathematically, the proposed DDUS, $\mathcal{U}$, is cast as follows:
\begin{equation}
\U \hspace{-0.05cm}=\hspace{-0.05cm} \Bigg\{\mathbf{U}\hspace{-0.05cm}\in\hspace{-0.07cm}\mathbb{R}^{\mathcal{|R|} \hspace{-0.05cm}\times\hspace{-0.05cm} \mathcal{|H|}} \hspace{-0.03cm}\Bigg{|} \mathbf{U} \hspace{-0.05cm}= \hspace{-0.04cm}\sum_{k = 1}^{K}\alpha_k \mathbf{U}_k, \sum_{k = 1}^{K}\alpha_k \hspace{-0.02cm}=\hspace{-0.02cm} 1,  \alpha_k \ge 0  \Bigg\}	\hspace{-0.09cm}\label{DDUS}\hspace{-0.09cm}
\end{equation}

\noindent where $\{\mathbf{U}_k\}_{k = 1,\ldots,K}$ is the set of daily REG profiles corresponding to $K$ previous days, which is hereinafter referred to as $\mathcal{S}$. Note that $\mathcal{U}$ is the convex hull of $\mathcal{S}$, thereby representing the smallest convex set that contains every scenario in $\mathcal{S}$. Moreover, by definition, the vertexes of the polyhedral uncertainty set $\mathcal{U}$ are in $\mathcal{S}$ \cite{bertsimas1997introduction}. Additionally, it is worth highlighting that the conservativeness level of the proposed DDUS is solely modeled by $K$. As a consequence, this parameter must be adjusted based on an out-of-sample test according to the decision maker's preference on cost and reliability.

\begin{figure}[h]
	\centering 	
	\includegraphics[width=7.3cm,height=3.5cm]{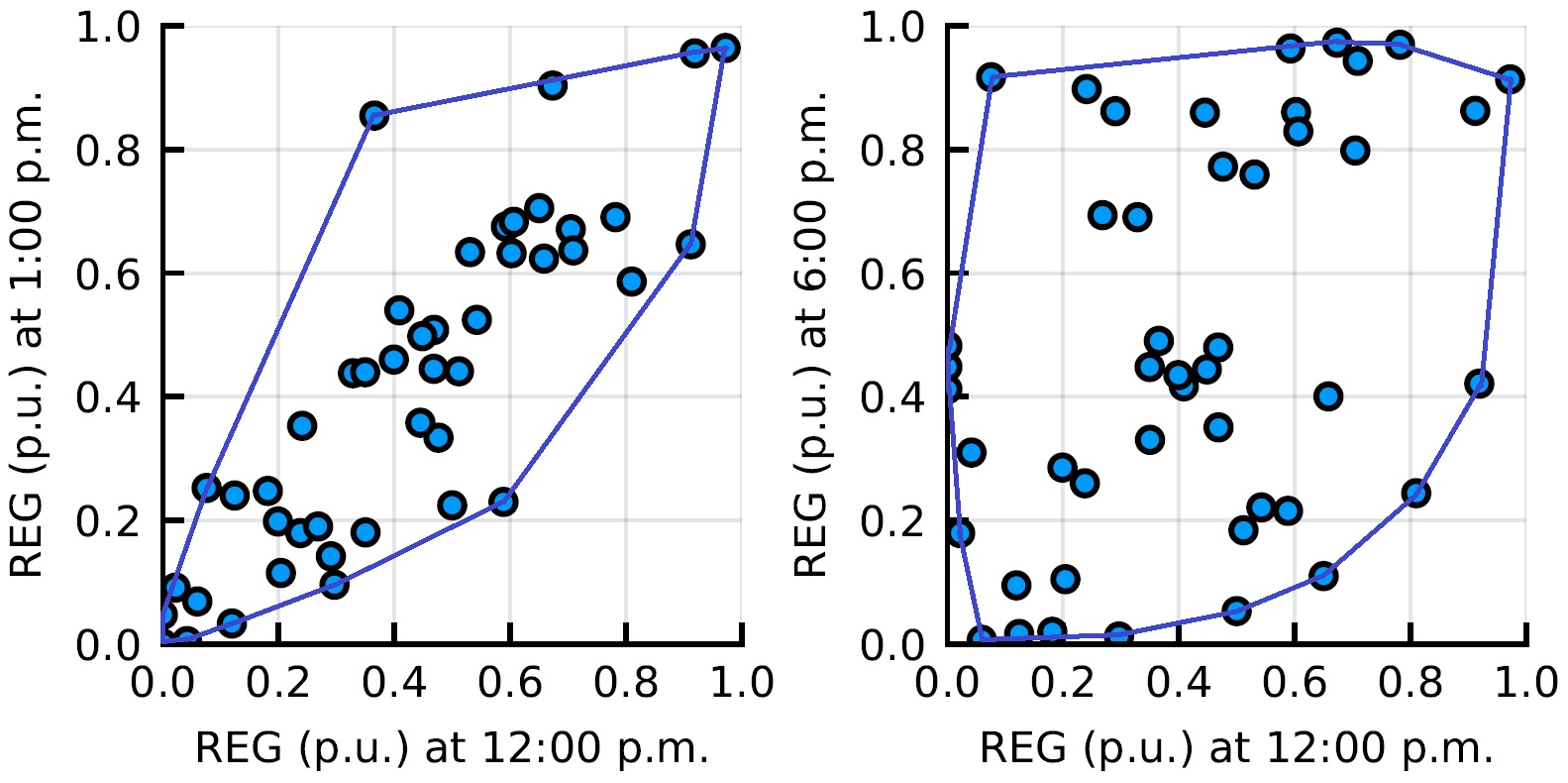}
	\vspace{-0.20cm}
	\caption{Examples of temporal dependences.}
	\label{fig:SameHour}
	\vspace{-0.1cm}
\end{figure}

\begin{figure}[h]
	\centering 	
	\includegraphics[width=7.3cm,height=3.5cm]{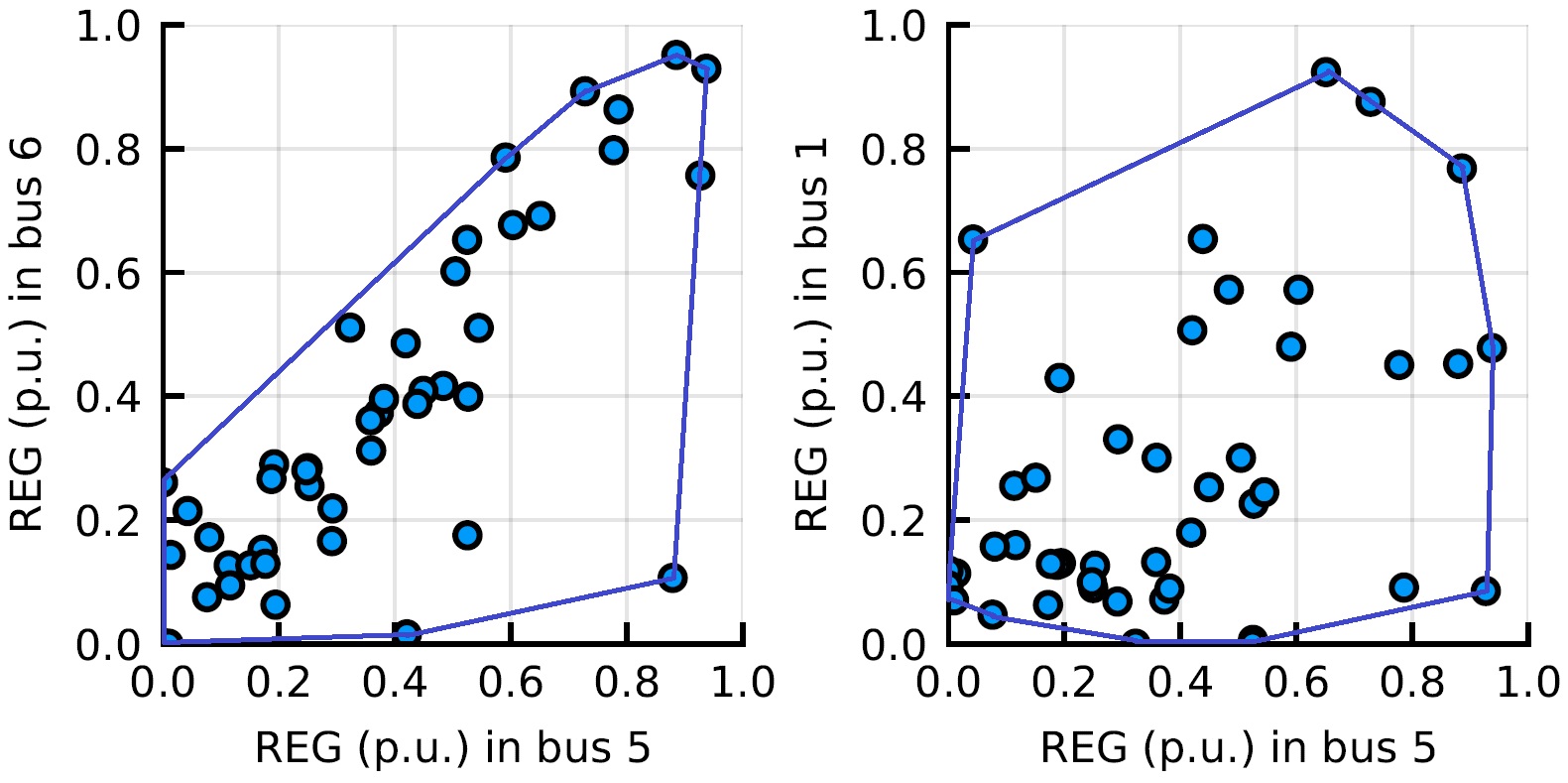}
	\vspace{-0.20cm}
	\caption{Examples of spatial dependences.}
	\label{fig:SameBus}
\end{figure}

The capability of DDUS to capture temporal and spatial dependences is illustrated in Figs. \ref{fig:SameHour} and \ref{fig:SameBus}, where two-dimensional projections of the multidimensional DDUS are shown for the data of the second case study examined in Section \ref{sec.Case}, for the period between 10/31/2012 and 12/15/2012. 
Fig. \ref{fig:SameHour} shows a high-dependence pattern, exhibited in bus 5 for 12:00 p.m. and 1:00 p.m., and a low-dependence pattern, observed in the same bus for 12:00 p.m. and 6:00 p.m. Likewise, different spatial dependences can be visualized in Fig. \ref{fig:SameBus} for 9:00 a.m.

An interesting interpretation for the use of observed generation profiles as vertexes of the uncertainty set is that we are implicitly performing an endogenous stress test for each feasible solution of the RUC problem whereby relevant multidimensional dependences found in the true hourly REG within a day are empirically accounted for. Furthermore, by using the proposed DDUS within an adaptive decision-making scheme, the evolution of the REG dependences across time can be  dynamically updated. Hence, strong assumptions about the nature of the true model behind data are not required. 

Moreover, note that the proposed framework can be adapted to consider sophisticated REG forecasts. As an example, any exogenous day-ahead forecasting procedure could be used to detrend/deseasonalize historical data and generate forecast errors. The errors computed for the $K$ previous days could thus be added to an updated day-ahead forecast to generate data-driven adjusted scenarios to build $\mathcal{S}$. With this simple preprocessing scheme, the shaping capability of the proposed data-driven approach could be straightforwardly combined with state-of-the-art forecasting methods, which will be addressed in future works.  

It is also worth mentioning that, if, instead of using real data, our data-driven model were built based on a Monte Carlo simulation of a time-series model, as done under a stochastic programming framework, our model can be viewed as a sampled version of a chance-constrained model. For instance, according to Calafiore \cite{Calafiore2009}, if the number of scenarios $K$ (generated through Monte Carlo sampling methods) in DDUS is greater than $|\mathcal{R}| |\mathcal{H}|/\varepsilon - 1$, where $\varepsilon \in (0,1)$ is a user-defined parameter that can be made arbitrarily small, 
 the expected probability of violations (shedding load or curtailing REG) remains below $\varepsilon$. 
As per our numerical experiments, the CCGA, in practice, converges with a much smaller subset of the $K$ sampled scenarios. Hence, in the case where our model is used to reproduce sample-based chance-constrained problems as done in \cite{Calafiore2009}, our robust approach would be a computationally efficient solution methodology.

\vspace{-0.20cm}
\section{Two-Stage Robust Unit Commitment Model}
\label{sec.Model}

For expository purposes, we use a standard two-stage robust optimization framework with a single uncertainty set \cite{ShabbirBook}, which is a particular instance of the general notion of stochastic robust optimization presented in \cite{BoZeng2015} and \cite{liu2016}. This model is simpler to describe and analyze yet bringing out the main features of the proposed data-driven approach. Based on previous works on robust energy and reserve scheduling \cite{Noemi2017,moreira2015energy,cobos2018network} and industry practice \cite{arroyo2005energy,galiana2005scheduling}, a base-case dispatch for energy is considered while ensuring $\epsilon$-feasibility for all realizations within the uncertainty set. The mathematical formulation for the proposed DDRUC problem is as follows:
\vspace{-0.05cm}
\begin{IEEEeqnarray}{ll}
	\min_{\mathcal{X},[\mathbf{\boldsymbol{\theta}}_{h},\mathbf{f}_{h}]_{h\in\mathcal{H}}}
	\,\sum _{h \in \mathcal{H}} \Big{\{} \mathbf{e}'(\mathbf{c}^{sd}_h + \mathbf{c}^{su}_h) + (\mathbf{c}^{g})' \mathbf{g}_h   \vspace{-0.5cm}\notag\\ \:\;\;\;\quad\qquad\qquad\qquad+ (\mathbf{c}^{up})' \mathbf{r}_h^{up}  + (\mathbf{c}^{dn})' \mathbf{r}_h^{dn} \Big{\}}      \label{eq.Master.Obj} \\
	 \text{subject to:} \notag \\
	\mathbf{A} \mathbf{f}_{h} + \mathbf{B} \mathbf{g}_{h} + \mathbf{P} \mathbf{\hat{u}}_{h} = \mathbf{{d}}_{h}  &\forall h \in \mathcal{H}	\label{eq.Master.BalancoEnergia}\\
	-\mathbf{\overline{f}} \leq \mathbf{f}_{h} \leq \mathbf{\overline{f}} &\forall h \in \mathcal{H} \label{eq.Master.PFLimit} \\
	\mathbf{f}_{h} = \mathbf{S}  \mathbf{\boldsymbol{\theta}}_h  \label{eq.Master.sec.Kirc} &\forall h \in \mathcal{H}\\
	\mathbf{\underline{G}} \mathbf{v}_h + \mathbf{r}_h^{dn} \le \mathbf{g}_{h} \leq \mathbf{\overline{G}} \mathbf{v}_h - \mathbf{r}_h^{up} &\forall h \in \mathcal{H} \label{eq.Master.En1}\\
	\mathbf{0} \leq \mathbf{r}_h^{up} \leq \mathbf{R}^{up} \: \mathbf{v}_h &\forall h \in \mathcal{H}\\
	\mathbf{0} \leq \mathbf{r}_h^{dn} \leq \mathbf{R}^{dn} \: \mathbf{v}_h &\forall h \in \mathcal{H} \label{eq.Master.En4}\\
	-\mathbf{b}(\mathbf{v}_h,\mathbf{v}_{h-1}) \le \mathbf{g}_{h}  - \mathbf{g}_{h-1}  \leq \mathbf{a}(\mathbf{v}_h,\mathbf{v}_{h-1}) &\forall h \in \mathcal{H} \label{Master.Ramp1}\\
	\{\mathbf{c}_{h}^{sd}, \mathbf{c}_{h}^{su},\mathbf{v}_{h}\}_{h\in\mathcal{H}} \in \mathcal{F} \quad \label{eq.Master.feasibilityset}\\
	\Phi(\mathcal{X}) \le \epsilon \label{tolerance} \\ 
	\Phi(\mathcal{X}) = \max_{\mathbf{U} \in \U} \Big\{ \min_{\substack{[\boldsymbol{\theta}^{wc}_h,\mathbf{f}^{wc}_h,\\\mathbf{g}^{wc}_h,\mathbf{s}^{+}_h,\\\mathbf{s}^{-}_h]_{h \in \mathcal{H}}}} \,\sum _{h \in \mathcal{H}}  \mathbf{e}'(\mathbf{s}_h^{+} + \mathbf{s}_h^{-}) \vspace{-0.1cm} \label{Oraculo}\\
	\quad\text{subject to:} \notag \\
	\quad\mathbf{A} \mathbf{f}_{h}^{wc} + \mathbf{B} \mathbf{g}_{h}^{wc} + \mathbf{P} \mathbf{u}_{h} = \notag \\ 
	\quad\qquad\qquad\qquad\mathbf{{d}}_{h} + \mathbf{s}_{h}^{+} - \mathbf{s}_{h}^{-}  \;\;\,: (\boldsymbol{\beta}_h) & \forall h \in \mathcal{H}	\label{eq.redispacho.BalancoEnergia}\\
	\quad-\mathbf{\overline{f}} \leq \mathbf{f}_{h}^{wc} \leq \mathbf{\overline{f}} \;\;\,\qquad\qquad\qquad: (\boldsymbol{\sigma}_h,\boldsymbol{\pi}_h) & \forall h \in \mathcal{H}\\
	\quad\mathbf{f}_{h}^{wc} = \mathbf{S} \mathbf{\boldsymbol{\theta}}_h^{wc} \,\,\,\quad\qquad\qquad\qquad: (\boldsymbol{\omega}_h)  &\forall h \in \mathcal{H}\\
	\quad\mathbf{g}_{h}^{wc} - \mathbf{g}_{h-1}^{wc}  \leq  \mathbf{a}(\mathbf{v}_h,\mathbf{v}_{h-1}) \quad\:: (\boldsymbol{\zeta}_h) &\forall h \in \mathcal{H} \\
	\quad\mathbf{g}_{h-1}^{wc} - \mathbf{g}_{h}^{wc} \leq  \mathbf{b}(\mathbf{v}_h,\mathbf{v}_{h-1}) \quad\,: (\boldsymbol{\kappa}_h) &\forall h \in \mathcal{H} \label{eq.redispacho.Rampa2}\\
	\quad\mathbf{g}_{h} - \mathbf{r}_{h}^{dn} \leq \; \mathbf{g}_{h}^{wc} \; \leq \mathbf{g}_{h} + \mathbf{r}_{h}^{up} \;: (\boldsymbol{\gamma}_h,\boldsymbol{\tau}_h) &\forall h \in \mathcal{H} \label{redispacho.limReserva}\\
	\quad\mathbf{0} \leq \mathbf{{s}}_{h}^{-} \leq \mathbf{d}_{h} \,\;\;\,\,\qquad\qquad\qquad:(\boldsymbol{\varsigma}_h) &\forall h \in \mathcal{H}\\
	\quad\mathbf{0} \leq \mathbf{{s}}_{h}^{+} \leq \mathbf{P} \mathbf{u}_{h} \;\qquad\qquad\qquad:(\boldsymbol{\xi}_h) &\forall h \in \mathcal{H} \Big\}\,\: \label{lastconstraint}
	\end{IEEEeqnarray}
\noindent where $\mathcal{X}$ is the set of first-level variables related to generation, i.e., $\mathcal{X} =[ \mathbf{c}^{sd}_h ,\mathbf{c}^{su}_h, \mathbf{g}_h, \mathbf{r}_h^{dn}, \mathbf{r}_h^{up}, \mathbf{v}_h ]_{h\in\mathcal{H}}$, and third-level dual variables are shown in parentheses. \label{setFirstStage}

\vspace{-0.04cm}
The two-stage DDRUC (\ref{eq.Master.Obj})--(\ref{lastconstraint}) is formulated as a (min-max-min) trilevel optimization problem. The first optimization level (\ref{eq.Master.Obj})--(\ref{tolerance}) determines the on/off statuses of generating units, as well as the energy and reserve scheduling. The objective function (\ref{eq.Master.Obj}) comprises shut-down costs, start-up costs, production costs, and up- and down-reserve costs. Constraint (\ref{eq.Master.BalancoEnergia}) represents nodal power balance under a dc power flow model. Expression (\ref{eq.Master.PFLimit}) represents transmission line power flow limits, while Kirchhoff's second law is accounted for through (\ref{eq.Master.sec.Kirc}). The limits for generation levels and up- and down-spinning reserves are imposed in expressions (\ref{eq.Master.En1})--(\ref{eq.Master.En4}). Inter-period ramping limits are modeled by (\ref{Master.Ramp1}), where the components of the auxiliary vector functions $\mathbf{a}(\cdot)$ and $\mathbf{b}(\cdot)$ are \cite{Noemi2017}: $a_{i}({v}_{ih},{v}_{i h-1})=RU_{i} {v}_{i h-1} + {SU_{i}} ({v}_{ih}-{v}_{i h-1}) + {\overline{G}_i}(1-{v}_{ih})$ and $b_{i}({v}_{ih},{v}_{i h-1}) = RD_{i} {v}_{ih}  + {SD_i}  ({v}_{ih-1}-{v}_{ih}) + {\overline{G}_i}(1-{v}_{ih-1})$. Following \cite{carrion2006computationally}, equation (\ref{eq.Master.feasibilityset}) represents constraints related to shut-down costs, start-up costs, and minimum up and down times. Equation (\ref{tolerance}) ensures redispatch capability within a feasibility tolerance, $\epsilon$, under the set $\mathcal{X}$ for all plausible realizations within the uncertainty set.

The second-level problem comprises the outer maximization problem in \eqref{Oraculo}. The goal of the second-level problem is to find the weights $\{\alpha_k\}_{k=1,...,K}$ in \eqref{DDUS} corresponding to the worst-case uncertainty realization $\mathbf{U} \in \U$ for a given value of $\mathcal{X}$. The measure of worst case is given by the minimum power imbalance function, which receives as inputs the values of the first- and second-level variables, $\mathcal{X}$ and $\mathbf{U}$, respectively. The third-level problem, i.e., the inner minimization problem (\ref{Oraculo})--(\ref{lastconstraint}), plays the role of the minimum power imbalance function
. The objective of this problem is to find a redispatch solution that minimizes the total sum of the mismatch variables artificially introduced in the power balance constraint (\ref{eq.redispacho.BalancoEnergia}). The mismatch variables, $\mathbf{s}^-_h$ and $\mathbf{s}^+_h$, can be interpreted as load shedding and REG spillage, respectively.
Constraints (\ref{eq.redispacho.BalancoEnergia})--(\ref{eq.redispacho.Rampa2}) are analogous to first-level expressions (\ref{eq.Master.BalancoEnergia})--(\ref{eq.Master.sec.Kirc}) and (\ref{Master.Ramp1}). The maximum operational deviation from nominal scheduled generation is controlled by constraint (\ref{redispacho.limReserva}), whereby generation redispatch levels $\mathbf{g}_{h}^{wc}$ are limited by up- and down-reserves. 

Within this approach, by ensuring $\epsilon$-feasibility for a given set of observed data points, the proposed DDRUC provides solutions that are optimized and tight for a given stress test setup. Therefore, the resulting reserve procurement can be interpreted as the least-cost schedule that allows reserve deliverability for realistic system-stress conditions. The proposed approach differs from current industry practice, where stress tests are applied \emph{ex post} as offline validation procedures. 

\vspace{-0.10cm}
\section{Solution Methodology}
\label{sec.Meth}
\vspace{-0.10cm}

This work leverages from the fact that the two lowermost levels \eqref{Oraculo}--\eqref{lastconstraint} correspond to a maximization, within a polyhedral uncertainty set $\U$, of a convex function given by the output of the inner minimization in (\ref{Oraculo}) as a function of $\mathbf{U}$. Therefore, from standard results of convex analysis, one of the vertexes belongs to the optimal solution set \cite{bertsimas1997introduction}. In other words, at the optimal solution, all $\alpha_k$ are equal to $0$ except that corresponding to the worst-case vertex, which is equal to $1$, thereby being binary valued. Hence, we can replace the continuous polyhedral uncertainty set $\U$ in the outer maximization of {(\ref{Oraculo})} with the discrete set of scenarios $\mathcal{S}$. By doing so, problem \eqref{eq.Master.Obj}--\eqref{lastconstraint} can be cast as a single-level MILP-based equivalent wherein expressions \eqref{Oraculo}--\eqref{lastconstraint} are replaced with one set of redispatch constraints \eqref{eq.redispacho.BalancoEnergia}--\eqref{lastconstraint} for each one of the $K$ scenarios in $\mathcal{S}$. Such a scenario-based model is structurally similar to that presented in \cite{Blanco2017} for a different problem, namely stochastic unit commitment. Unfortunately, addressing such a full-scenario-based equivalent by the branch-and-cut algorithm may lead to intractability.

Alternatively, problem \eqref{eq.Master.Obj}--\eqref{lastconstraint}, and its single-level equivalent, likewise, are suitable for the CCGA \cite{BoZeng2011}, which ensures finite convergence to global optimality by iteratively solving a master problem and an oracle subproblem. At each iteration, the master problem, which is a relaxed version of the full-scenario-dependent equivalent problem, finds a trial solution that is evaluated in terms of operational feasibility by the oracle procedure. The oracle returns to the master problem violated constraints that will change the master problem solution at the next iteration. The algorithm terminates when the oracle does not find any violated constraint, thereby certifying the incumbent trial solution as globally optimal.

\vspace{-0.3cm}
\subsection{Master Problem}

The master problem is a relaxed version of problem \eqref{eq.Master.Obj}--\eqref{lastconstraint}. Thus, at each iteration $j$ of the CCGA, the master problem is an instance of MILP as follows:
\vspace{-0.10cm}
\begin{IEEEeqnarray}{ll}
	\min_{\substack{\mathcal{X},[\mathbf{\boldsymbol{\theta}}_{h},\mathbf{f}_{h}]_{h\in\mathcal{H}}\\     [   \boldsymbol{\theta}^{(k)}_h,\mathbf{f}^{(k)}_h,\mathbf{g}^{(k)}_h,\mathbf{s}^{+(k)}_h,\\\mathbf{s}^{-(k)}_h]_{h\in\mathcal{H}, k=1,\ldots,|\mathcal{S}_j|}}}
	\;\sum _{h \in \mathcal{H}} \Big{\{} \mathbf{e}'(\mathbf{c}^{sd}_h + \mathbf{c}^{su}_h) + (\mathbf{c}^{g})' \mathbf{g}_h  \;\;\;\; \vspace{-0.5cm}\notag\\ \;\;\;\qquad\qquad\qquad\qquad\;\;\;\;\;\;\;+ (\mathbf{c}^{up})' \mathbf{r}_h^{up}  + (\mathbf{c}^{dn})' \mathbf{r}_h^{dn} \Big{\}} \;\;     \label{SD.Obj} \qquad \quad  \vspace{-0.1cm}\\
	\text{subject to:} \notag \\
	\;\text{Constraints (\ref{eq.Master.BalancoEnergia})--(\ref{eq.Master.feasibilityset})}\label{SD.1}\\
	\; \text{Redispatch constraints for scenario $k,\; k=1,\ldots,|\mathcal{S}_j|$} \label{SD.redispatch}\\
	\;\sum _{h \in \mathcal{H}}  \mathbf{e}'(\mathbf{s}_h^{+(k)} + \mathbf{s}_h^{-(k)}) \le \epsilon, \quad  k=1,\ldots,|\mathcal{S}_j|  \label{SD.Oraculo}
\end{IEEEeqnarray}

\noindent where set $\mathcal{S}_j$ comprises the scenarios of $\mathcal{S}$ selected by the oracle subproblem until iteration $j$ of the CCGA. 

Expressions \eqref{SD.Obj}--\eqref{SD.1} are identical to \eqref{eq.Master.Obj}--\eqref{eq.Master.feasibilityset}. As per \eqref{SD.redispatch}, a set of redispatch constraints is iteratively added. Such redispatch constraints correspond to \eqref{eq.redispacho.BalancoEnergia}--\eqref{lastconstraint} where 
variables $\boldsymbol{\theta}^{wc}_h,\mathbf{f}^{wc}_h,\mathbf{g}^{wc}_h,\mathbf{s}^{+}_h,$ and $\mathbf{s}^{-}_h$ are respectively replaced with new variables $\boldsymbol{\theta}^{(k)}_h,\mathbf{f}^{(k)}_h,\mathbf{g}^{(k)}_h,\mathbf{s}^{+(k)}_h,$ and $\mathbf{s}^{-(k)}_h$, whereas $\mathbf{u}_h$ is replaced with the corresponding scenario of $\mathcal{S}_j$ selected at previous iterations by the oracle subproblem. Finally, in \eqref{SD.Oraculo}, the system power imbalance under every scenario in $\mathcal{S}_j$ is bounded by the threshold $\epsilon$.

\vspace{-0.20cm}
\subsection{Oracle Subproblem}
\label{Oraclesubsection}

The oracle subproblem corresponds to the two lowermost optimization levels (\ref{Oraculo})--(\ref{lastconstraint}) for given upper-level decisions provided by the preceding master problem. This bilevel program can be solved in two different ways:

\begin{enumerate}

\item By the enumeration of all feasible values for variables $\alpha_k$. This solution approach, hereinafter referred to as the \emph{inspection-based oracle}, consists in the serial or parallel evaluation of the imbalance through the solution of the third-level problem for all $K$ scenarios in $\mathcal{S}$. The inspection-based oracle is thus equivalent to searching through all vectors $\{\alpha_k\}_{k = 1,\ldots,K}$ that are candidates for optimality, i.e., those with one entry equal to $1$ and all other entries equal to $0$. Note that each optimization problem is a linear program, which runs in polynomial time with interior point methods. As the value of $K$ is not related to the instance size, the inspection-based oracle can be easily implemented as a computational routine that receives a pair $(\mathcal{X}^*,\mathcal{S})$ and returns the worst-case scenario $\mathbf{U}^*$ in $\mathcal{S}$ and the associated power imbalance $\Phi(\mathcal{X}^*)$ in polynomial time \cite{Bertsimas2017}.
	
\item By the application of branch-and-cut algorithms (readily available in off-the-shelf MILP solvers) to a single-level MILP-based equivalent subproblem. Based on a discrete representation of $\mathcal{S}$, this solution approach consists in casting the original bilevel subproblem as a single-level MILP equivalent, denoted by \emph{MILP-based oracle} subproblem. Such an equivalent requires modeling variables $\alpha_k$ as binary, using the dual of the lower level of the subproblem, and applying well-known integer algebra results \cite{floudas1995} to recast the resulting bilinear terms as linear expressions. The MILP-based oracle subproblem is thus a linearized version (relying on standard disjunctive constraints) of this mixed-integer nonlinear problem:
\vspace{-0.20cm}
\end{enumerate}
	\begin{IEEEeqnarray}{ll}
		\Phi(\mathcal{X}) =\max_{\substack{[\alpha_k]_{k\in\mathcal{K}},\\ [\boldsymbol{\beta}_h,\boldsymbol{\gamma}_h,\boldsymbol{\zeta}_h,\\\boldsymbol{\kappa}_h,\boldsymbol{\xi}_h,\\\boldsymbol{\pi}_h,\boldsymbol{\sigma}_h,\\\boldsymbol{\varsigma}_h,\boldsymbol{\tau}_h,\\\boldsymbol{\omega}_h,\mathbf{u}_h]_{h\in\mathcal{H}}}} \sum _{h \in \mathcal{H}} \Big{\{} \boldsymbol{\beta}_h'  (\mathbf{{d}}_{h} - \mathbf{P} \mathbf{u}_{h})  
		\vspace{-1.5cm} \notag\\
		\;\,\,\qquad\qquad\qquad\qquad\qquad - \boldsymbol{\pi}_h' \mathbf{\overline{f}} - \boldsymbol{\sigma}_h' \mathbf{\overline{f}}  \notag\\
		\;\,\,\qquad\qquad\qquad\qquad\qquad+\boldsymbol{\gamma}_h' (\mathbf{g}_{h} - \mathbf{r}_{h}^{dn}) \notag\\
		\;\,\,\qquad\qquad\qquad\qquad\qquad-\boldsymbol{\tau}_h'  (\mathbf{g}_{h} + \mathbf{r}_{h}^{up})\notag\\
		\;\,\,\qquad\qquad\qquad\qquad\qquad- \boldsymbol{\zeta}_h'  \mathbf{a}(\mathbf{v}_h,\mathbf{v}_{h-1}) \; \notag\\
		\;\,\,\qquad\qquad\qquad\qquad\qquad+ \; \boldsymbol{\kappa}_h'  \mathbf{b}(\mathbf{v}_h,\mathbf{v}_{h-1})
		\notag\\
		\;\,\,\qquad\qquad\qquad\qquad\qquad-\boldsymbol{\varsigma}_h' \mathbf{{d}}_{h} \;   -\;\boldsymbol{\xi}_h'\mathbf{P} \mathbf{u}_{h} \Big{\}}\;\; \label{DEDproblem}\vspace{-0.2cm}\\
		\text{subject to:}  \notag \\
		\mathbf{u}_h = \sum_{k \in \mathcal{S}}\alpha_k \mathbf{u}_{hk} \;\; \qquad\qquad\qquad\qquad \forall h\in\mathcal{H} \label{SOS1a} \\
		\sum_{k \in \mathcal{S}}\alpha_k = 1 \label{SOS1b}\\
		\alpha_k \in \{0,1\} \;\qquad\qquad\qquad\qquad\qquad\:\, \forall  k\in\mathcal{S} \label{SOS1c}  \\
		\boldsymbol{\beta}_h - \boldsymbol{\varsigma}_h \le \mathbf{e} \,\,\;\;\;\qquad\qquad\quad\qquad\qquad \forall h \in \mathcal{H} \label{dual1}\\
		-\boldsymbol{\beta}_h - \boldsymbol{\xi}_h \le \mathbf{e} \qquad\qquad\qquad\quad\qquad\,\, \forall h \in \mathcal{H}\\
		- \mathbf{A}' \boldsymbol{\beta}_h + \boldsymbol{\pi}_h - \boldsymbol{\sigma}_h
		- \boldsymbol{\omega}_h = \mathbf{0} \quad \;\;\;\,\quad \forall h \in \mathcal{H}\\
		- \mathbf{S}' \boldsymbol{\omega}_h = \mathbf{0} \,\quad\qquad\,\;\quad \qquad\qquad\qquad\forall h \in \mathcal{H}\\
		-\boldsymbol{\gamma}_h+\boldsymbol{\tau}_h-\mathbf{B}' \boldsymbol{\beta}_h + \boldsymbol{\zeta}_h  \notag\\\quad \qquad- \boldsymbol{\zeta}_{h+1} - \boldsymbol{\kappa}_h +\boldsymbol{\kappa}_{h+1}    \geq \mathbf{0} \quad\:\,\,\,\, \forall h \in \mathcal{H}\\
		\boldsymbol{\pi}_h , \boldsymbol{\sigma}_h,\boldsymbol{\tau}_h,\boldsymbol{\gamma}_h,\boldsymbol{\zeta}_h, \boldsymbol{\kappa}_h, \boldsymbol{\varsigma}_h, \boldsymbol{\xi}_h \geq \;\mathbf{0} \quad \forall h \in \mathcal{H}. \label{DEDproblemEnd}
	\end{IEEEeqnarray}
	
	In \eqref{DEDproblem}--\eqref{DEDproblemEnd}, the objective function \eqref{DEDproblem} and constraints \eqref{dual1}--\eqref{DEDproblemEnd} represent, respectively, the objective function and constraints of the dual formulation for the third-level problem presented in (\ref{Oraculo})--(\ref{lastconstraint}). The proposed DDUS is equivalently modeled by constraints \eqref{SOS1a}--\eqref{SOS1c}, where $\mathbf{u}_{hk}$ stands for the $h$-th column of the $k$-th data point, $\mathbf{U}_{k}$. According to \eqref{DDUS} and the aforementioned results of convex analysis, $\mathcal{U}$ is replaced with $\mathcal{S}$ in \eqref{SOS1a} and \eqref{SOS1b}, whereas variables $\alpha_k$ are characterized as binary in \eqref{SOS1c}. Note that problem \eqref{DEDproblem}--\eqref{DEDproblemEnd} is a  mixed-integer nonlinear optimization problem with bilinear terms in the objective function \eqref{DEDproblem}. Such nonlinearities involve products of $\u_h$ and third-level dual variables. As is customary in the literature \cite{Noemi2017,moreira2015energy, cobos2018network}, an equivalent MILP formulation for \eqref{DEDproblem}--\eqref{DEDproblemEnd}, i.e., the MILP-based oracle, is achieved by applying standard integer algebra results \cite{floudas1995} that allow recasting the nonlinear products in \eqref{DEDproblem} as linear expressions.

\vspace{-0.20cm}
\subsection{Algorithm}
\label{subsec.CCGA}

For given $\mathcal{S}$ and $\epsilon$, the proposed CCGA works as follows:
\begin{algorithm}[H]
\renewcommand{\thealgorithm}{}
{\small \caption{{\small CCGA($\mathcal{S}$, $\epsilon$})} \label{Alg.CCGA}
	\begin{algorithmic}[1]
		\State{Initialization: $j\leftarrow0$ and $\mathcal{S}_j\leftarrow \emptyset$.}
		\State{Solve the master problem \eqref{SD.Obj}--\eqref{SD.Oraculo} to obtain $\mathcal{X}^* $.}
		\State{Solve the subproblem for ($\mathcal{X}^*$, $\mathcal{S}$)} to obtain $\Phi(\mathcal{X}^*)$ and the worst-case scenario $\mathbf{U}^*$.
		\If{$\Phi(\mathcal{X}^*) \le \epsilon$}\textbf{:} STOP 
		\Else{}\textbf{:} $j\leftarrow j+1$, $\mathcal{S}_j\leftarrow \mathcal{S}_{j-1} \cup \{\mathbf{U}^*\}$, and go to step 2
		\EndIf.
	\end{algorithmic}
} 
\end{algorithm}
\vspace{-0.20cm}
The algorithm iteratively adds violated constraints into the master problem and terminates when infeasibility is within $\epsilon$.

In order to reduce the search space, vertex-identification algorithms \cite{RikoJacob2002} could be implemented \emph{ex ante} to determine the subset of $\mathcal{S}$ comprising the vertexes of $\mathcal{U}$. However, in a high-dimensional case, it is very unlikely to find a point in $\mathcal{S}$ that is not a vertex of $\mathcal{U}$. 
Hence, the use of vertex-identification procedures, in practice, may deteriorate the performance of the CCGA as only an insignificant number of points of $\mathcal{S}$, most likely none, would be excluded for being in the interior of the DDUS.

\vspace{-0.25cm}
\section{Numerical Results}
\label{sec.Case}

This section reports results from an illustrative 4-bus system and the IEEE 118-bus test system over a 24-hour time span. For the sake of reproducibility, data for both test systems can be downloaded from \cite{DataDrivenSystemLAMPS}. For expository purposes, wind-related uncertainty is considered. The source of data for wind power generation is the Global Energy Competition (GEFCom) \cite{hong2014global,hong2016probabilistic}. 

As explained in the next section, the proposed DDRUC has been assessed with two benchmark models based on two-stage robust optimization and two-stage stochastic programming. Such an assessment has been conducted through backtesting on historical data observed over a given set of days. 

The benchmark RUC model, hereinafter referred to as BRUC, relies on the budget-constrained uncertainty set formulated in Section \ref{sub.Conventional}. DDRUC has been implemented for different time windows, i.e., for different values of $K$, whereas different combinations of budgets $\Gamma$ and $\Lambda$ were considered for BRUC. For quick reference, the instances of both models are denoted by DDRUC($K$) and BRUC($\Gamma$,$\Lambda$), respectively. For a given day $d$, the uncertainty set for DDRUC($K$) is directly defined using the REG profiles observed in the $K$ previous days, i.e., considering $\mathcal{S}$ as $\{\mathbf{U}_k\}_{k=d-K,\ldots,d-1}$. 

For the sake of a fair comparison, 
BRUC is also adaptively adjusted across the rolling-horizon study, similarly to DDRUC, in order to prevent over-conservative solutions. This is done through a moving window of previous observed days within which the lower and upper box-like limits used in the budget-constrained uncertainty sets are defined according to the corresponding hourly minimum and maximum production limits for each renewable unit. Our tests indicated that the length of such a moving window is not as relevant for BRUC as the selection of $K$ is for DDRUC. Hence, based on trial and error, the budget-constrained uncertainty sets for BRUC($\Gamma$,$\Lambda$) were built using the previous $35$ days. 

The instances of DDRUC were addressed by the CCGA with both the inspection- and the MILP-based oracle subproblems. In contrast, BRUC was solved by the CCGA with a modified version of the MILP-based oracle subproblem, as done in \cite{Noemi2017}. To that end, expressions \eqref{SOS1a} and \eqref{SOS1b} in problem \eqref{DEDproblem}--\eqref{DEDproblemEnd}, representing the uncertainty set constraints, were replaced with \eqref{BUSlin1}--\eqref{BUSlin3}, whereas the linearization scheme mentioned in Section \ref{Oraclesubsection} was applied to the resulting bilinear terms. 

Additionally, in order to provide a comparison between the proposed data-driven robust approach and stochastic programming, a benchmark stochastic unit commitment model, hereinafter referred to as BSUC, was analyzed. BSUC is a mixed-integer linear program structurally similar to the DDRUC master problem \eqref{SD.Obj}--\eqref{SD.Oraculo}. Three differences characterize BSUC, namely 1) an expected imbalance cost term is added to the objective function \eqref{SD.Obj}, 2) the set of scenarios in \eqref{SD.redispatch} comprises uncertainty realizations obtained from a scenario-generation procedure using probabilistic information, and 3) imbalance requirements are eliminated by dropping constraint \eqref{SD.Oraculo}. The value of the imbalance cost was set to \$500/MWh, which corresponds approximately to 25 times the average value (or 13 times the highest value) of the fuel costs, $\mathbf{c}^{g}$, for the thermal generators.

The scenario-generation procedure for BSUC was based on a multivariate lognormal distribution inferred from the previous 100 days. First, the historical data were normalized and a logarithmic transformation was applied. Then, a multivariate normal distribution was estimated and 500 independent and equally distributed scenarios were generated based on the expected value vector and the covariance matrix. Finally, the inverse transformations were applied to rescale the scenarios. As for the solution methodology, BSUC was tackled by the branch-and-cut algorithm, which was initialized with the previous day's solution for binary scheduling variables and reserve contributions. 

Simulations were run using Gurobi 7.0.2 under JuMP (Julia 0.5) on a Xeon E5-2680 processor at 2.5 GHz and 128 GB of RAM.

\vspace{-0.3cm}

\subsection{Evaluation Methodology}
\label{sub.BenchmarkingModel}

The performances of DDRUC($K$), BRUC($\Gamma$,$\Lambda$), and BSUC have been compared in terms of their tradeoff between cost and robustness. The evaluation methodology consists in conducting a rolling-horizon out-of-sample backtest over a set of days for which realistic REG data are available. 

For each day $d$ within the backtest horizon, the instances under assessment are solved using available historical REG data for past days, i.e., excluding day $d$. The solution of such instances of DDRUC, BRUC, and BSUC yields the corresponding values of the total cost for that day. Subsequently, solution robustness is quantified by computing the infeasibility of the generation dispatch associated with the resulting generation schedules and the actual REG scenario observed for day $d$. To that end, based on industry practice, a single-period version of the inner minimization problem in \eqref{Oraculo}--\eqref{lastconstraint} is solved for each hour of this day. Hence, based on the resulting hourly levels of load shedding and REG spillage, out-of-sample statistics for robustness are devised.
\vspace{-0.20cm}
\subsection{Illustrative 4-Bus System}
\label{subsec.4bus}

\begin{figure}[t]
\centering 	
\includegraphics{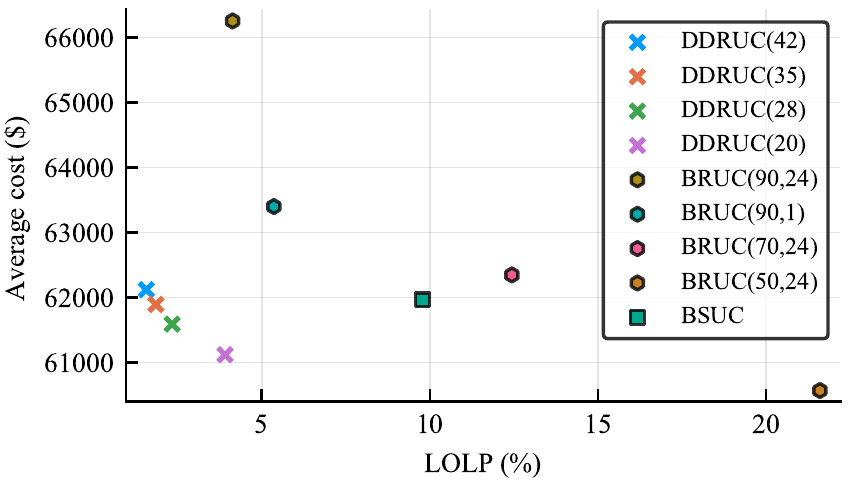}
\vspace{-0.20cm}
\caption{Average cost vs. LOLP for the 4-bus system.}
\label{fig:4busParetoKgleast}
	\vspace{-0.15cm}
\end{figure}

First, we consider an illustrative test system consisting of $4$ buses, $4$ transmission lines, $14$ thermal generators, and $2$ wind farms \cite{DataDrivenSystemLAMPS}. Wind generation data were rescaled from buses $2$ and $5$ from \cite{hong2014global}. The optimality gap of Gurobi was set at $0\%$ and the CCGA was implemented until no imbalance was found for the instances of DDRUC and BRUC. Hence, such simulations were run to optimality for $\epsilon = 0$. As for BSUC, a time limit of 1,500 seconds was considered.
\begin{table}[t]
\centering
\caption{4-Bus System --Results from DDRUC, BRUC, and BSUC}	  	
\vspace{-0.150cm}
\small{
	\begin{tabular}{ p{1.67cm} p{1.5cm} p{0.68cm} p{0.81cm} p{1.5cm} p{1.5cm} }
		\toprule\noalign{\vskip 0.45mm}
		&\multicolumn{1}{c}{Avg.}&\multicolumn{1}{c}{Avg.}&\multicolumn{1}{c}{Avg.} &\multicolumn{2}{|c}{$ $}\\
		\multicolumn{1}{c}{Instance}&\multicolumn{1}{|c}{time}&\multicolumn{1}{c}{iter.}& \multicolumn{1}{c}{cost}&\multicolumn{1}{|c}{LOLP}&\multicolumn{1}{c}{PWS}\\
		& \multicolumn{1}{|c}{(s)}&\multicolumn{1}{c}{no.} &\multicolumn{1}{c}{(\$)} &\multicolumn{1}{|c}{(\%)}&\multicolumn{1}{c}{(\%)}\\
		\hline\noalign{\vskip 0.45mm}
		\multicolumn{1}{l}{DDRUC(42)}&\multicolumn{1}{|r}{$20.0\,(34.7)$} &\multicolumn{1}{r}{$13.4$}&$62122$&\multicolumn{1}{|r}{$1.56$}&\multicolumn{1}{r}{$2.67$}\\
		DDRUC(35)&\multicolumn{1}{|r}{$17.7\,(28.8)$} &\multicolumn{1}{r}{$12.8$} &$61894$&\multicolumn{1}{|r}{$1.85$}&\multicolumn{1}{r}{$3.53$}\\
		DDRUC(28)&\multicolumn{1}{|r}{$15.2\,(23.2)$} &\multicolumn{1}{r}{$12.0$}& $61592$&\multicolumn{1}{|r}{$2.33$} &\multicolumn{1}{r}{$4.15$}\\
		DDRUC(20)&\multicolumn{1}{|r}{$11.9\,(17.0)$} & \multicolumn{1}{r}{$10.9$}& $61121$&\multicolumn{1}{|r}{$3.91$} &\multicolumn{1}{r}{$5.43$}\\
		BRUC(90,24)&\multicolumn{1}{|r}{$157.7\qquad\;\;\:$}&\multicolumn{1}{r}{$6.2$}&$66256$&\multicolumn{1}{|r}{$4.13$}&\multicolumn{1}{r}{$1.00$}\\
		BRUC(90,1) &\multicolumn{1}{|r}{$688.7\qquad\;\;\:$}&\multicolumn{1}{r}{$71.5$}&$63401$&\multicolumn{1}{|r}{$5.36$}&\multicolumn{1}{r}{$1.09$}\\
		BRUC(70,24)&\multicolumn{1}{|r}{$115.9\qquad\;\;\:$}&\multicolumn{1}{r}{$6.5$}&$62349$&\multicolumn{1}{|r}{$12.44$}&\multicolumn{1}{r}{$2.86$}\\
		BRUC(50,24)&\multicolumn{1}{|r}{$13.0\qquad\;\;\:$}& \multicolumn{1}{r}{$6.2$}&$60571$&\multicolumn{1}{|r}{$21.60$}&\multicolumn{1}{r}{$5.25$}\\
		BSUC&\multicolumn{1}{|r}{$1494.0\qquad\;\;\:$}& \multicolumn{1}{r}{$-$}&$61967$&\multicolumn{1}{|r}{$9.78$}&\multicolumn{1}{r}{$2.96$}\\
				\bottomrule
	\end{tabular}
	\vspace{-0.15cm}
	\label{Tabela1}}
\end{table}
The evaluation backtest was conducted for $336$ days of realistic wind power generation data. Table \ref{Tabela1} and Fig. \ref{fig:4busParetoKgleast} summarize the results from the backtest for several instances of DDRUC, BRUC, and BSUC. 

Columns $2$--$4$ of Table \ref{Tabela1} respectively provide the average computing times, iteration numbers, and costs over the $336$ runs of the models being examined. For the instances of DDRUC, the first figure in column $2$ corresponds to the method using the inspection-based oracle subproblem, whereas the second figure in parentheses is associated with the approach relying on the solution of the MILP-based equivalent for the oracle subproblem. Column $5$ presents a reliability index referred to as loss of load probability (LOLP), which is defined as the fraction of hours with load shedding exceeding $0.1\%$ of the system load. Analogously, column $6$ reports the probability of wind spillage (PWS), defined as the fraction of hours in which wind spillage exceeds $0.1\%$ of the available wind power generation. Thus, the values of LOLP and PWS listed in Table \ref{Tabela1} represent valuable statistical information on the robustness of the solutions provided by each model.

Table \ref{Tabela1} shows the impact of $K$ on the computational effort required by DDRUC. It should be emphasized that the number of iterations of the CCGA grows at a less than linear rate with $K$ whereas computing times behave similarly. As a relevant result, the impact of $K$ on computing times is lower when the subproblem is solved by the vertex inspection, which outperforms the MILP solver for all instances.

As compared with BRUC, DDRUC is in general more effective from a computational perspective. For the instances with no temporal budget constraints ($\Lambda = 24$), the average numbers of iterations for BRUC are approximately half of those required to solve the instances of DDRUC. However, such reduced numbers of iterations do not imply shorter computing times because the MILP-based oracle for BRUC is computationally expensive. Moreover, for $\Lambda = 1$, which is far less conservative, BRUC requires significantly more iterations to converge on average, thereby increasing the computational burden. An exception is the instance BRUC(50, 24), which is a particular case with no temporal budget constraints and a spatial budget corresponding to precisely a single renewable generator deviating from the expected production. Interestingly, this was the case for which BRUC achieved the worst solution robustness in terms of LOLP, which is 5.5 times higher than the worst level attained by DDRUC. Regarding BSUC, its required computational effort is considerably greater than those of DDRUC and BRUC. The 1,500-s time limit was reached almost every day of the rolling-window study.      

As can also be observed in Table \ref{Tabela1}, as $K$ grows from 20 to 42, both LOLP and PWS decrease. Note also that, in comparison with BRUC and BSUC, all solutions provided by DDRUC gave rise to acceptably small LOLP values ranging between $1.56$\% and $3.91$\%, while keeping the average cost at reasonable levels within a narrow $1.6$\% band. The superiority of DDRUC over BRUC and BSUC is thus evidenced in terms of the tradeoff between cost and robustness. Such a tradeoff is further illustrated in Fig. \ref{fig:4busParetoKgleast}, where robustness is expressed in terms of the out-of-sample reliability index LOLP. Note that the solutions to DDRUC dominate in Pareto sense (lower cost and lower LOLP) almost all solutions provided by BRUC. The only non-dominated BRUC instance is the aforementioned BRUC(50, 24), which, albeit incurring the lowest average cost, leads to an unreasonably high LOLP. Analogously, despite featuring a reasonable average cost, the solution provided by BSUC gave rise to a very high value of LOLP, namely 9.78\%, which was substantially outperformed by all instances of DDRUC. Hence, the BSUC solution was dominated by the instances DDRUC(20), DDRUC(28), and DDRUC(35).     

\begin{table}[!t]
\centering
\caption{4-Bus System -- Average Cost Breakdown ($10^3$ \$)}	  
\vspace{-0.20cm}	
\small{
	\begin{tabular}{l|cccc}
		\toprule
		\multirow{ 2}{*}{\hspace{0.10cm} Instance}& \multirow{ 2}{*}{Total}&\multirow{ 2}{*}{ Production} & \multirow{2}{*}{Reserve}&	\multicolumn{1}{c}{Start-up and}\\
		&  &  &  & \multicolumn{1}{c}{shut-down}\\
		\midrule
		DDRUC(35)    &61.9&55.7&6.1&0.1\\
		DDRUC(20)    &61.1&55.5&5.5&0.1\\
		BRUC(90,24) &66.3&57.2&9.0&0.1\\
		BRUC(90,1)  &63.4&56.2&7.1&0.1\\
		BSUC  		 &62.0&55.7&6.2&0.1\\
		\toprule
	\end{tabular}%
	\label{Tabela2}}
\vspace{-0.15cm}
\end{table}

\begin{figure}[t]
	\centering 
	\includegraphics{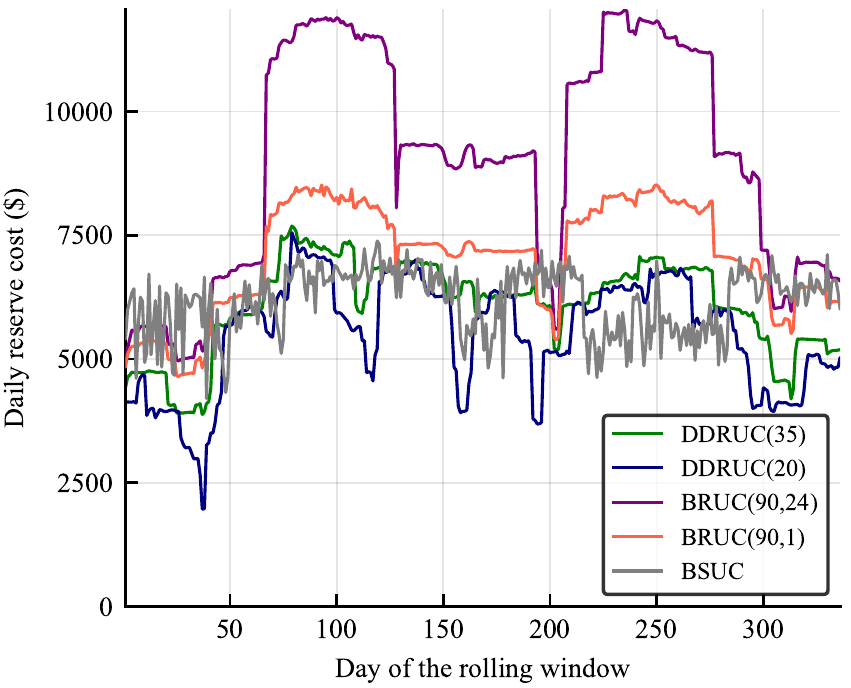}
	\vspace{-0.20cm}
	\caption{Daily reserve costs for the 4-bus system.}
	\label{fig:Reserve}
		\vspace{-0.25cm}
\end{figure}

A breakdown of the average costs provided in column 4 of Table \ref{Tabela1} is presented in Table \ref{Tabela2} for representative instances of DDRUC, BRUC, and BSUC models. Thus, average production costs, average reserve costs, and average start-up and shut-down costs are reported. It can be observed that the largest cost component is the average production cost, representing between 86\% and 91\% of the average total cost for all instances. However, the average reserve cost, which amounts to between 9\% and 14\% of the average total cost, is the most relevant in terms of cost difference for any pair of instances. Overall, reserve costs for DDRUC and BSUC are considerably lower than those for BRUC. For further details, daily reserve costs are shown in Fig. \ref{fig:Reserve}. As compared to BRUC, DDRUC and BSUC were able to address the next-day REG profile with lower levels of reserves. Moreover, daily reserve costs for DDRUC are less variable than those for BSUC due to the inherent stochastic nature of the scenario-generation procedure.  

\begin{figure}[t!]
\centering 
\includegraphics{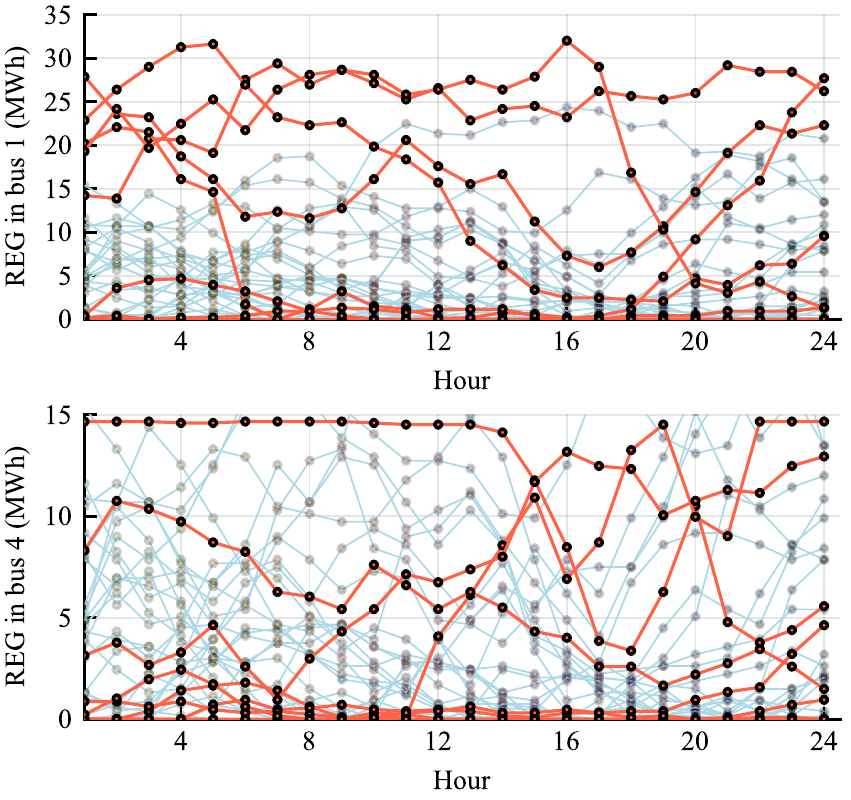}
\vspace{-0.20cm}
\caption{Scenarios in $\mathcal{S}$ for the $43^{rd}$ day of the rolling window for the 4-bus system.}
\label{fig:Scenarios}
	\vspace{-0.20cm}
\end{figure}

\begin{figure}[t]
\centering 
\includegraphics{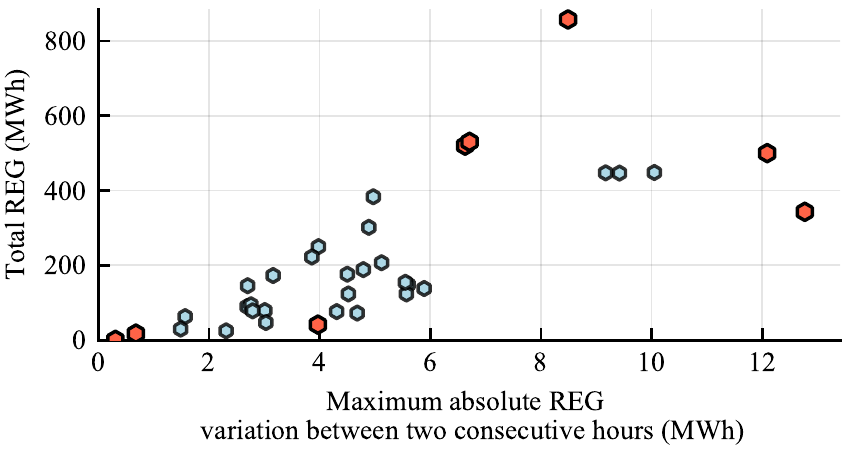}
\vspace{-0.20cm}
\caption{Scatter plot for the scenarios selected by DDRUC(35) in the $43^{rd}$ day of the rolling window for the 4-bus system.}
\label{fig:Scenarios2}
	\vspace{-0.20cm}
\end{figure}

For illustration purposes, Fig. \ref{fig:Scenarios} depicts two two-dimensional projections of the scenarios in $\mathcal{S}$ for the $43^{rd}$ day of the rolling window. The 8 scenarios colored in red correspond to the scenarios selected by the CCGA for DDRUC(35). Not surprisingly, some of the selected scenarios present very high or very low REG profiles. Such scenarios are more likely to yield higher imbalance costs (from either wind spillage or load shedding) at intermediate iterations of the CCGA. Salient features of the selected scenarios are also characterized by the scatter plot in Fig. \ref{fig:Scenarios2}, which shows the total 24-hour REG and the largest absolute variation of renewable generation between two consecutive hours of each scenario. Note that the two scenarios with the largest absolute REG variation between two consecutive periods were also selected. These scenarios also represent a challenge to the system as they stress its ramping capability and require higher reserve levels. 

\subsection{118-Bus System}
\label{subsec.118bus}
The second case study is a modified version of the IEEE $118$-bus test system with $10$ wind farms. The optimality gap of Gurobi was set at $0.05\%$ for the master problem, the subproblem was solved to optimality through the vertex inspection, and $\epsilon$ was equal to $0.1\%$ of the total system demand. 

For this case study, we first assessed the performance of DDRUC($35$) with a second instance, denoted by DDRUC($35$+$10$), whereby intra-day 
dynamics of wind power generation for the same period in the previous year were considered. Thus, the scenarios for DDRUC($35$+$10$) were associated with the previous $35$ days and a ten-day period from the year before, beginning five days prior to the date corresponding to the day under analysis in the backtest scheme. For both instances, the backtest was performed for $214$ days under realistic wind power generation data \cite{hong2016probabilistic}. The results listed in Table \ref{Tabela118bus} indicate that the proposed DDRUC can provide a robust and cost-effective generation schedule. Moreover, it is interesting to note that the cost increase due to the consideration of $10$ additional vertexes from the previous year is negligible, while the values of LOLP and PWS respectively drop from $3.2$\% and $5.8$\% to $2.7$\% and $4.6$\%.

\begin{table}[t]
\centering
\caption{118-Bus System -- Results from the 214-Day Backtest}
\small
\begin{tabular}{ p{2.73cm} p{0.68cm} p{0.68cm} p{0.1cm}p{0.68cm} p{0.68cm} }
	\toprule
	&\multicolumn{2}{r}{$\;\;\;$DDRUC(35)$\;\;\;$}&&\multicolumn{2}{r}{DDRUC(35+10)}\\
	\cline{2-6}\noalign{\vskip 0.25mm}
	Average time (s) & \multicolumn{2}{r}{$774\qquad$}  & & \multicolumn{2}{r}{$792\qquad$} \\
	Maximum time (s) & \multicolumn{2}{r}{$1752\qquad$}  & & \multicolumn{2}{r}{$2220\qquad$} \\
	Average no. of iter. & \multicolumn{2}{r}{$8.5\qquad$}  & & \multicolumn{2}{r}{$8.4\qquad$}\\
	Average cost ($10^4\,$\$) &\multicolumn{2}{r}{$110.4\qquad$}& &\multicolumn{2}{r}{$110.5\qquad$} \\
	LOLP (\%)              &\multicolumn{2}{r}{$3.2\qquad$} & & \multicolumn{2}{r}{$2.7\qquad$}\\
	PWS (\%)              &\multicolumn{2}{r}{$5.8\qquad$} & & \multicolumn{2}{r}{$4.6\qquad$}\\
	\bottomrule
\end{tabular}
\small
\label{Tabela118bus}
\end{table}


\begin{table}[t]
	\centering
	\caption{118-Bus System -- Results from the 10-Day Backtest}	  	
	\small{
		\begin{tabular}{ p{1.8cm} p{0.68cm} p{0.85cm} p{0.39cm} p{0.185cm} p{0.185cm} }
			\toprule\noalign{\vskip 0.45mm}
			&\multicolumn{1}{|c}{Avg.}&\multicolumn{1}{c}{Avg.}&\multicolumn{1}{c}{Avg.} &\multicolumn{2}{|c}{$\,$}\\
			\multicolumn{1}{c}{Instance}&\multicolumn{1}{|c}{time}&\multicolumn{1}{c}{iter.}& \multicolumn{1}{c}{cost}&\multicolumn{1}{|r}{LOLP}&PWS\\
			& \multicolumn{1}{|c}{(s)}&\multicolumn{1}{c}{no.} &\multicolumn{1}{c}{$\,$($10^4\,$\$)$\,$} &\multicolumn{1}{|c}{$\;\;$(\%)}&\multicolumn{1}{c}{(\%)}\\
			\hline\noalign{\vskip 0.45mm}
			\multicolumn{1}{l}{DDRUC(42)}&\multicolumn{1}{|r}{$494$} &\multicolumn{1}{c}{$5.0$}&\multicolumn{1}{c}{$110.2$}&\multicolumn{1}{|r}{$1.67$}&\multicolumn{1}{c}{$0.83$}\\
			\multicolumn{1}{l}{DDRUC(35+10)}&\multicolumn{1}{|r}{$696$} &\multicolumn{1}{c}{$5.7$}&\multicolumn{1}{c}{$110.2$}&\multicolumn{1}{|r}{$1.67$}&\multicolumn{1}{c}{$0.83$}\\
			DDRUC(35)&\multicolumn{1}{|r}{$378$} &\multicolumn{1}{c}{$5.6$} &\multicolumn{1}{c}{$110.1$}&\multicolumn{1}{|r}{$1.67$}&\multicolumn{1}{r}{$1.25$}\\
			DDRUC(28)&\multicolumn{1}{|r}{$804$} &\multicolumn{1}{c}{$6.6$}&\multicolumn{1}{c}{$109.6$}&\multicolumn{1}{|r}{$1.25$} &\multicolumn{1}{r}{$5.00$}\\
			BRUC(70,24)&\multicolumn{1}{|r}{$3847$}&\multicolumn{1}{c}{$3.0$}&\multicolumn{1}{c}{$109.6$}&\multicolumn{1}{|r}{$8.33$}&\multicolumn{1}{r}{$1.67$}\\
			BRUC(50,24)&\multicolumn{1}{|r}{$3741$}& \multicolumn{1}{c}{$3.0$}&\multicolumn{1}{c}{$108.3$}&\multicolumn{1}{|r}{$19.17$}&\multicolumn{1}{r}{$5.83$}\\
			\bottomrule
			\vspace{-0.1cm}
		\end{tabular}
		\label{Tabela3}}
\end{table}

It is also relevant to highlight the low computational burden exhibited by both DDRUC instances. The significance of such a result was validated by performing the 214-day backtest with BRUC, which was the only benchmark model used for this case study due to the poor computational behavior featured by BSUC for the illustrative example. Unfortunately, for BRUC, the CCGA failed to converge in practical computing times.

DDRUC has been further assessed with BRUC through a reduced backtest relying on a randomly selected sequence of $10$ days. Simulations were run for a feasibility tolerance $\epsilon$ of $0.1$\% of the system demand and a practical 1-hour time limit. For the instances of BRUC examined, the time limit was reached during the solution of the master problem. For such instances, the CCGA was allowed to run until the corresponding iteration was completed. Table \ref{Tabela3} summarizes the results from the 10-day backtest. The reported results suggest that the proposed DDRUC is a scalable alternative for RUC.

\vspace{-0.2cm}

\section{Conclusion}
\label{sec.Conclusion}

This paper has described a new, comprehensive, and parsimonious method to characterize REG uncertainty in day-ahead RUC. The proposed method relies on an alternative scenario-based polyhedral uncertainty set that is built through a novel data-driven approach. Unlike conventional budget-constrained uncertainty sets, the proposed polyhedral uncertainty set is characterized directly from data through the convex hull of a set of previously observed REG profiles. Thus, relevant empirical information regarding the existing complex and time-varying dynamics of REG sources is embedded in the vertexes of the polyhedron. In addition, a salient feature of the proposed uncertainty set is its capability to consider any exogenous model (or scheme) for adjusting, including or excluding (filtering or clustering) scenarios based on existing methods used in industry. Moreover, the resulting robust counterpart for generation scheduling can be efficiently solved by the CCGA until $\epsilon$-global optimality. This relevant practical aspect stems from the reduced complexity of the oracle subproblem, which is solvable in polynomial time. For the case studies analyzed in the paper, a relatively small number of past renewable generation profiles were required by the proposed data-driven approach to outperform two benchmarks, namely an RUC model based on a conventional budget-constrained uncertainty set and a stochastic unit commitment model. Numerical results indicate that the proposed method might be a practical, scalable, easy to specify, and cost-efficient alternative tool for managing wind power variability in the unit commitment problem.

Future research will consider 1) exogenous day-ahead forecasting methods for nominal scenarios, and 2) data-processing schemes based on forecasting errors to generate the vertexes for the uncertainty set. In addition, alternative budget-constrained uncertainty sets will be analyzed. Another relevant avenue of research is the extension of the proposed data-driven approach within a stochastic robust optimization framework.

\vspace{-0.15cm}

\bibliographystyle{IEEEtran}
\bibliography{bibliography}

\begin{IEEEbiography}
	[{\includegraphics[width=1in,height=1.25in,clip]{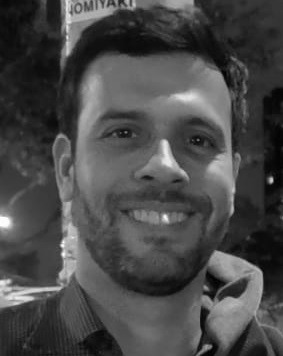}}]
	{Alexandre Velloso} (S'17) received his M.Sc. degree in mathematical methods in finance from the National Institute for Pure and Applied Mathematics (IMPA), Rio de Janeiro, Brazil in 2005. Since 2010, he has been working for the Brazilian Innovation Agency -- Financiadora de Estudos e Projetos (Finep). He is pursuing a Ph.D. degree in electrical engineering at the Pontifical Catholic University, Rio de Janeiro, Brazil. Currently, he is a visiting scholar at the School of Industrial and Systems Engineering at the Georgia Institute of Technology, Atlanta, USA. 
\end{IEEEbiography}

\begin{IEEEbiography}
	[{\includegraphics[width=1in,height=1.25in,clip,keepaspectratio]{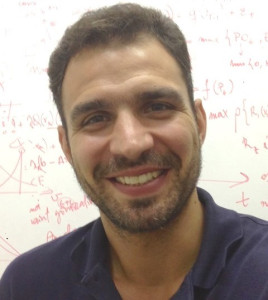}}]
	{Alexandre Street} (S'06--–M'10--–SM'17) holds an M.Sc. and a D.Sc. (PhD) in electrical engineering (emphasis in operations research) from the Pontifical Catholic University of Rio de Janeiro (PUC-Rio), Rio de Janeiro, Brazil. From 2003 to 2007, he participated in several projects related to strategic bidding in the Brazilian energy auctions and market regulation at the Power System Research Consulting (PSR), Rio de Janeiro. From August 2006 to March 2007, he was a visiting researcher at the Universidad de Castilla-La Mancha, Ciudad Real, Spain. At the beginning of 2008, he joined the Electrical Engineering Department at PUC-Rio, where he is currently an Associate Professor and teaches energy economics, optimization and statistics. He is one of the founders of the Laboratory of Applied Mathematical Programming and Statistics (LAMPS), where he is currently the research director. Alexandre Street is a senior member of the IEEE Power and Energy Society and CNPq senior scientist.
\end{IEEEbiography}

\begin{IEEEbiography}
		[{\includegraphics[width=1in,height=1.25in,clip]{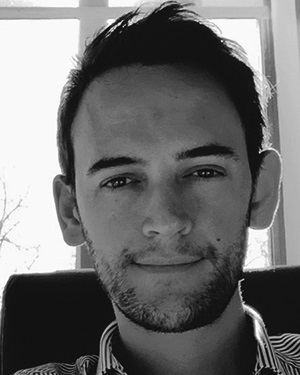}}]
		{David Pozo} (S'06--–M'13--–SM'18) received his B.S. and Ph.D. degrees in electrical engineering from the University of Castilla-La Mancha, Ciudad Real, Spain, in 2006 and 2013, respectively. Since 2017, he is Assistant professor at the Skolkovo Institute of Science and Technology (Skoltech), Moscow, Russia. Prior to Skoltech Dr. Pozo worked as a postdoctoral fellow at the Pontifical Catholic University of Chile and the Pontifical Catholic University of Rio de Janeiro. 
		His research interest lies in the field of power systems and includes operations research, uncertainty, game theory, and electricity markets. He also focuses on problems of optimization and flexibility of modern power systems. 
		Since 2018, Prof. Pozo is leading the research group on Power Markets Analytics, Computer Science and Optimization (PACO). 
\end{IEEEbiography}

\begin{IEEEbiography}
		[{\includegraphics[width=1in,height=1.25in,clip, keepaspectratio]{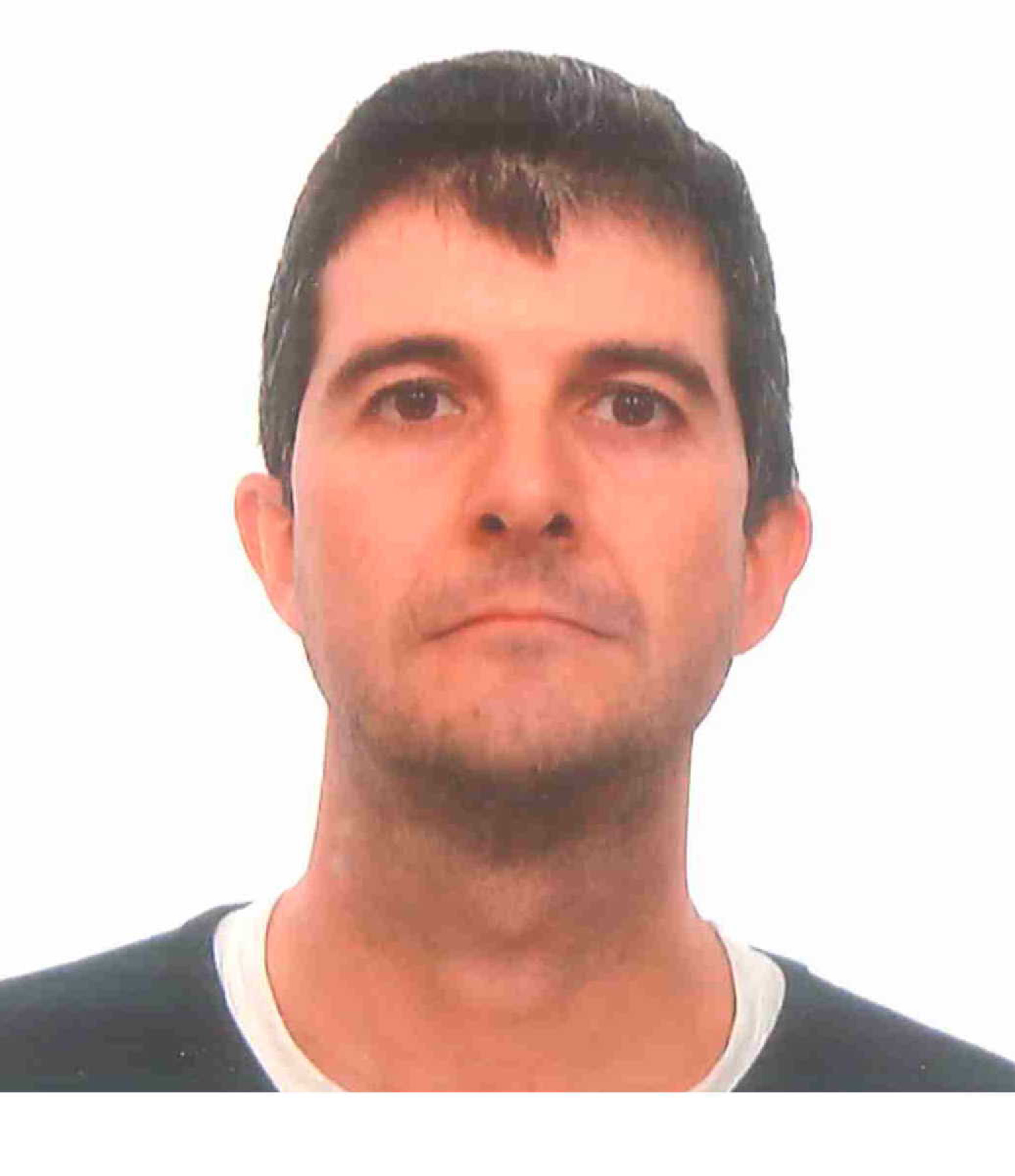}}]
		{Jos\'e M. Arroyo} (S'96--–M'01--–SM'06) received the Ingeniero Industrial degree from the Universidad de M\'alaga, M\'alaga, Spain, in 1995, and the Ph.D. degree in power systems operations planning from the Universidad de Castilla-La Mancha, Ciudad Real, Spain, in 2000.
	
	From June 2003 to July 2004, he held a Richard H. Tomlinson Postdoctoral Fellowship with the Department of Electrical and Computer Engineering,
	McGill University, Montreal, QC, Canada.
	From January to August 2018, he held a visiting professorship with the same institution.
	
	He is currently a Full Professor of electrical engineering with the Universidad de Castilla-La Mancha. His research interests include operations, planning, and economics of power systems, as well as optimization.
\end{IEEEbiography}

\begin{IEEEbiography}
	[{\includegraphics[width=1in,height=1.25in,clip]{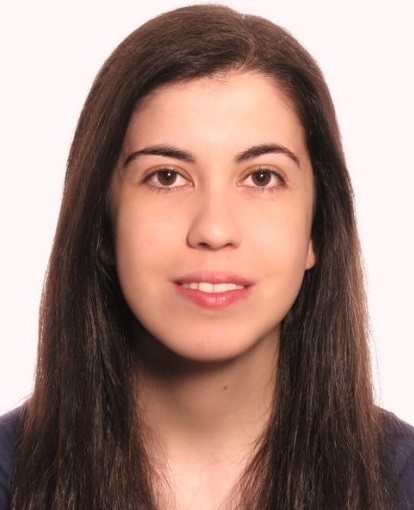}}]
	{Noemi G. Cobos}(S'15) received the Ingeniero Industrial, M.Sc., and Ph.D. degrees in electrical engineering from the Universidad de Castilla-La Mancha, Ciudad Real, Spain, in 2013, 2014, and 2018, respectively. She is currently working with the SmartGrids department at the Instituto Tecnol\'{o}gico de la Energ\'{i}a (Valencia, Spain). Her research interests are in the fields of power system operation and economics.
\end{IEEEbiography}

\end{document}